\documentclass{svjour3}
\usepackage{amscd,amsfonts,amssymb,amsmath,latexsym,array,hhline,xcolor,graphicx}

\usepackage{latexsym}
\usepackage{amssymb}

\newcommand\F{\mbox{I\kern-2pt F}}
\newcommand\cA{{\cal A}}
\newcommand\cE{{\cal E}}
\newcommand\cC{{\cal C}}
\newcommand\cF{{\cal F}}

\newcommand\cI{{\cal I}}

\newcommand\cL{{\cal L}}

\newcommand\cX{{\cal X}}

\newcommand\cO{{\cal O}}

\newcommand\cS{{\cal S}}

\newcommand\cY{{\cal Y}}

\newcommand\e{{\varepsilon}}
\newtheorem{theo}{Theorem}[section]
\newtheorem{prop}[theo]{Proposition}
\newtheorem{lemm}[theo]{Lemma}

\newcommand\fdem{$\Box$}
\newcommand\beq{\begin{equation}}
\newcommand\eeq{\end{equation}}
\newcommand\bea{\begin{eqnarray}}
\newcommand\eea{\end{eqnarray}}
\newcommand\bean{\begin{eqnarray*}}
\newcommand\eean{\end{eqnarray*}}

\begin{document}

\title{Consumption-Investment Problem with Transaction Costs
for L\'evy-Driven Price Processes }

\titlerunning{Consumption-Investment Problem with Transaction Costs}

\author{Dimitri De Valli\`ere \and 
Yuri Kabanov \and  Emmanuel L\'epinette }


\institute{
    D. De Valli\`ere \at
             Laboratoire de Math\'ematiques, Universit\'e de
Franche-Comt\'e, 16 Route de Gray, 25030 Besan\c{c}on, cedex,
France
 \and
 Yu. Kabanov  \at
              Laboratoire de Math\'ematiques, Universit\'e de
Franche-Comt\'e, 16 Route de Gray, 25030 Besan\c{c}on, cedex,
France,  and\\
National Research University  Higher School of Economics, International Laboratory of Quantitative Finance, Moscow, Russia. \\
              \email{Youri.Kabanov@univ-fcomte.fr}        
\and
E. L\'epinette \at Ceremade, Universit\'e Paris Dauphine, Place du Mar\'echal De Lattre De Tassigny, \\
75775 Paris cedex 16, France, and\\
National Research University  Higher School of Economics, International Laboratory of Quantitative Finance, Moscow, Russia.\\
\email Emmanuel.Denis@ceremade.dauphine.fr
}

\date{Received: date / Accepted: date}

\maketitle

\begin{abstract} We consider an optimal control problem for a linear stochastic integro-diffe\-rential equation with conic constraints on the phase variable and the control of singular-regular type. Our setting includes
consumption-investment problems for models of financial markets in the presence of proportional
transaction costs where  the price of the assets are given by   a geometric L\'evy
process and  the investor is allowed to take short positions. We prove that the Bellman function of the problem is
a viscosity solution of the HJB equation. A uniqueness theorem for the solution of the latter is established.
Special attention is paid to the Dynamic Programming Principle.


\end{abstract}

\section{Introduction } In this paper we study  the classical consumption-investment model with infinite horizon in the presence of transaction costs. Our aim is to extend the results of \cite{K-Kl} to the case where the price evolution is given by a  geometric L\'evy process. Namely, we show that the Bellman function is a viscosity solution of the corresponding Hamilton--Jacobi--Bellman equation. We also prove a uniqueness theorem for the latter.

Mathematically, the consumption-investment problem with transaction costs  is a regular-singular control problem for a linear stochastic equation in a cone. Its specificity is that the Bellman function may not  be  smooth and, therefore, one cannot use the verification theorem
(at least, in its traditional form) because the It\^o formula cannot be applied. Nevertheless, one can show that the Bellman function is a solution of the HJB equation in viscosity sense. Though the general line of arguments  is common, one needs to re-examine each step of the proof. In particular, for the considered jump-diffusion model, the HJB equation contains an integro-differential operator and the test functions involved in the definition of the viscosity solution must be ``globally" defined. It seems that already in 1986 H.M. Soner noticed that the control problems with jump parts can
be considered in the framework of the theory of viscosity solutions,
\cite{Soner86a}, \cite{Soner86b}.

There is a growing literature on extension of the concept of viscosity solutions to  equations with integro-differential operators, see, e.g., \cite{Sayah},  \cite{Alv-Tourin}, \cite{Pham}, \cite{Barles-Imbert}, \cite{Barles-Chasseigne-Imbert}, \cite{Arisawa}, \cite{Arisawa2}.
There are  several variants of the definition of viscosity solution. Our choice is intended to serve the model with a positive utility function. The definition can be viewed as a simplified version of that adopted in \cite{JK}.

A rather detailed study of  consumption-investment problems
under transaction costs when the  prices follow exponential
L\'evy processes and the investor is constrained to keep long positions in all assets, money included, was undertaken in papers by Benth et {\it al.} \cite{Benthal} and \cite{Bentha}.
Our geometric approach seems to be more general than that
of the mentioned papers  where the authors consider a ''parametric" version of the stock market, with transactions always passing through money (i.e. either ``buy stock" or ``sell stock").
A more important difference is that in our setting the investor may take short positions as was always assumed in the classical papers
\cite{Merton}, \cite{DN}, \cite{Sh-Son}. If  short positions are admitted, the ruin may happen due to a jump of the price process.
That is why the natural, "classical", setting  considered  here leads to a different HJB equation of a more complicated structure.
Following the ideas from the paper \cite{K-Kl} we derive the Dynamic Programming Principle splitted into two separate assertions.  Though it is the principal tool which allows to check that the Bellman function is a viscosity solution of the HJB equation, it is rarely discussed
in the literature (and even taken for granted, see, e.g., in \cite{AMS}, \cite{Sh-Son}, \cite{Benthal}).

The main results of the paper is Theorem \ref{conditional} claiming that if the Bellman function is continuous up to the boundary then it is a viscosity
solution of the HJB equation and the
uniqueness theorem for the Dirichlet problem arising in the model, Th. \ref{unique}. We formulate the latter in terms of the Lyapunov function, an object that is defined in terms of the truncated operator, in which the utility function is not involved. Its introduction  allows us to disconnect the problems of the uniqueness of a solution and  the existence of a classical supersolution. 

Probably, the most important result of the paper is the uniqueness theorem for the HJB equation with a non-local operator.  
In contrast to the methods developed in \cite{BI} which are based on (very technical) extensions of the Ishii lemma     
we use the latter in its original and (very transparent) formulation. 

The structure of the problem is the following. In Sections 2 and 3 we introduce the model dynamics and describe the goal functional providing comments on the concavity of the Bellman function $W$.
In Section 4 we show that  the Bellman function, if finite, is continuous in the interior of the solvency cone. In Section 5  we give a formal description of the HJB equation. Sections 6 and 7 contain
a short account of basic facts on viscosity solutions for
integro-differential operators. In Section 8 we explain the role
of classical supersolutions to the HJB equations. Section 9 is devoted to the Dynamic Programming Principle. In Section 10 we use it to show that the Bellman function is the solution of our HJB equation. Section 11 contains a uniqueness theorem formulated in terms of a Lyapunov function. In Section 12 we provide examples of Lyapunov functions and classical supersolutions.


\section{The Model} Our setting  is
more general than that of the standard model of financial market
under constant proportional transaction costs. In particular, the
cone $K$ is not supposed to be polyhedral. We assume that the
asset prices  are geometric L\'evy processes.
 Our framework appeals to a  theory of
viscosity solutions for non-local integro-differential operators.

Let $Y=(Y_t)$ be an ${\bf R}^d$-valued semimartingale  on
 a stochastic basis
 $(\Omega, \cF,{\bf F},P)$  with the
  trivial initial $\sigma$-algebra. Let $K$ and $\cC$ be
 {\bf proper} closed cones
in ${\bf R}^d$ 
such that 
$\cC\subseteq {\rm int}\,K\neq\emptyset$.
Define the set $\cA$ of controls $\pi=(B,C)$ as the set of {\bf predictable} ${\bf R}^d$-valued
c\`adl\`ag processes of bounded variation such that,
up to an evanescent set, \beq \label{RN-constr}
  \dot B \in -K, \qquad  \dot C \in \cC.
\eeq
Here $\dot B$ denotes a (measurable version) of  the Radon--Nikodym derivative of $B$ with respect to the total variation process $||B||$. The notation $\dot C$ has a similar sense. Though models with arbitrary $C$
is of interest, we restrict ourselves in the present paper by considering consumption processes admitting intensity. To this end, we define
$\cA_a$ as the set of controls $\pi$ with absolutely continuous components $C$ such that
 the increment $C_0=0$. For the elements of $\cA_a$ we have
$c:=dC/dt\in \cC$.

The controlled process $V=V^{x,\pi}$ is the solution of
  the linear  system
 \beq
\label{dynam} dV^i_t=V^i_{t-} dY^i_t+dB^i_t-dC^i_t, \ \quad
V^i_{0-}=x^i,\quad i=1,...,d.
 \eeq
In general, $\Delta V_{0}=\Delta B_0$ is not  equal to zero: the investor may revise the portfolio when entering the market at time zero.

The solution of (\ref{dynam}) can be expressed explicitly
using the Dol\'eans-Dade exponentials
 \beq
\label{doleans}
\cE_t(Y^i)=e^{Y^i_t-(1/2)\langle Y^{ic}\rangle_t}\prod_{s\le t}(1+\Delta Y^i_s)e^{-\Delta Y^i_s}.
\eeq
Namely,
 \beq
\label{dynam1} V^i_t=\cE_t(Y^i)x^i+\cE_t(Y^i)
\int_{[0,t]}\cE_{s-}^{-1}(Y^i)(dB^i_s-dC^i_s), \quad i=1,...,d.
 \eeq

We introduce the stopping time
 $$
 \theta=\theta^{x,\pi}:=\inf\{t:\  V^{x,\pi}_t\notin {\rm int}\, K\}.
$$
For  $x\in {\rm int}\,K$ we consider the subsets $\cA^x$ and
$\cA^x_a$ of ``admissible" controls for which $\pi=I_{[0,\theta^{x,\pi}]}\pi$  and $\{V_-+\Delta B\in {\rm int}\, K\}=\{V_-\in {\rm int}\, K\}$. In financial context, $\theta$ is the time of ruin.  When $V^{x,\pi}$  leaves
the interior of the solvency cone the control of the portfolio and the  consumption stops. 
The process $V$ given by (\ref{dynam}), continues to evolve after the time $\theta$ but for us the relevance has only the stopped process
 $V^{x,\pi, \theta}$. 

It is natural to assume that the process $V$ does not leave the interior of $K$ due to a jump of $B$: the investor is reasonable enough not to ruin himself by making a too expensive portfolio revision.

\smallskip

The important hypothesis that the cone $K$ is proper, i.e. $K\cap
(-K)=\{0\}$, or equivalently, ${\rm int}\, K^*\neq \emptyset$,
corresponds
 to the  model of financial market with {\it efficient friction}. In a financial
context $K$ (usually containing ${\bf R}^d_+$) is interpreted as
the solvency region and $C=(C_t)$ as the consumption process; the
process $B=(B_t)$ describes accumulated  fund transfers. In the ``standard" model with proportional transaction costs (sometimes referred to as the model of currency market)
$$
K={\rm cone}\,\{(1+\lambda^{ij})e_i-e_j,\ e_i,\ 1\le i,j\le
d\}
$$
where $\lambda^{ij}\ge 0$ are transaction costs coefficients, see  Section 3.1 in the book \cite{K-Saf} for details and other examples.
\smallskip

The process $Y$ represents the relative price movements. If $S^i$ is the price process of the $i$th asset, then $dS^i_t=S^i_{t-}dY^i_t$ and $S^i_{t}=S^i_{0}\cE_t(Y^i)$.
Without loss of generality we  assume that  $S^i_{0}=1$ for all $i$. In this case $Y^i$ is the so-called stochastic logarithm of $S^i$.  The formula (\ref{dynam1}) can be re-written as follows:
 \beq
\label{dynam2} V^i_t=S^i_tx^i+S^i_t
\int_{[0,t]}\frac 1{S^i_{s-}}(dB^i_s-dC^i_s), \quad i=1,...,d.
 \eeq

We shall work  assuming that
\beq
\label{dY}
Y_t=\mu t + \Xi  w_t +\int_0^t\int z( p(dz,dt)-q(dz,dt))
\eeq
where $\mu\in {\bf R}^d$, $w$ is a $m$-dimensional standard Wiener
process and  $p(dz,dt)$ is a Poisson random measure with the compensator $q(dz,dt)=\Pi(dz)dt$ such that $\Pi(dz)$ is a measure concentrated on  $]-1,\infty[^d$. Note that the latter property of the L\'evy measure corresponds to the financially meaningful case where $S^i>0$.
For the $m\times d$-dimensional matrix $\Xi$ we put $A=\Xi\Xi^*$.
We assume that
\beq\label{Pi}
\int (|z|^2\wedge |z|)\Pi(dz)<\infty
\eeq
and this  assumption validates the formula (\ref{dY}): by definition,  
\bean
\int_0^t\int z( p(dz,dt)-q(dz,dt))&:= &\int_0^t\int_{\{|z|\le 1\}} z( p(dz,dt)-q(dz,dt))\\
&&+\int_0^t
\int
_{\{|z|> 1\}} z p(dz,dt)-\int_0^t \int_{\{|z|> 1\}}z\int q(dz,dt))
\eean
where the first integral is defined as a stochastic one, while the second and the third are 
the usual Lebesque integrals, both finite (a.s.). 

\smallskip

{\bf Notation.}  By typographical reasons we shall use the notation 
$D_x$ instead of common ${\rm diag }\,x$ for the diagonal operator (or matrix) generated by the vector $x=(x^1,...,x^d)$, i.e. 
$$
D_xz=(x^1z^d,...,x^dz^d). 
$$ 

The system (\ref{dynam}) can be written in the integral vector form as follows: 
\beq
\label{dynamvec} 
V_t=x+\int_0^t D_{V_{s-}} (\mu ds+\Xi dW_s) +\int_0^t\int D_{V_{s-}}z\left(p(dz,ds)-q(dz,ds)\right)+B_t- C_t.
\eeq

\smallskip 
It is important to note that the jumps of $Y$ and $B$ cannot occur
simultaneously. More precisely, the process $|\Delta B||\Delta Y|$ is indistinguishable of zero. Indeed, for any $\e>0$ we have, using the predictability of the process $\Delta B=B-B_-$,  that
\bean
E\sum_{s\ge 0}|\Delta B_s||\Delta Y_s|I_{\{|\Delta Y_s|>\e\}}&=&E\int_0^\infty\int|\Delta B_s|I_{\{|z|>\e\}}|z|p(dz,ds)\\
&=&E\int_0^\infty\int|\Delta B_s||z|I_{\{| z|>\e\}}\Pi (dz)ds=0
\eean
 because for each $\omega$ the set $\{s: \Delta B_s(\omega)\neq 0\}$ is at most countable and its Lebesgue measure is equal to zero.
 Thus, the process $|\Delta B||\Delta Y|I_{\{|\Delta Y|>\e\}}$ is indistinguishable of zero and so is the process $|\Delta B||\Delta Y|$.


\smallskip
   It follows that $\Delta B_\theta=0$. Since the predictable process
 $I_{\{V_-\in \partial K\}}I_{[0,\theta]}$ has at most countable number of jumps, the same reasoning as above leads to the conclusion that $I_{\{V_-\in \partial K\}}|\Delta Y|
I_{[0,\theta]}$ is  indistinguishable of zero. This means that $\theta$ is the first
 moment when either $V$ or $V_-$ leaves  ${\rm int}\,K$. This property will be used in the proof that $W$ is lower semicontinuous on ${\rm int}\, K$.

\smallskip
In our proof of the Dynamic Programming Principle (needed to
derive the HJB equation) we shall assume that the stochastic basis
is a canonical one, that is the space of c\`adl\`ag functions and $P$ is a measure 
under which the coordinate mapping is
  the L\'evy process.



\section{Goal Functionals and Concavity of the Bellman Function}

Let $U:\cC\to {\bf R}_+$ be a concave  function such that $U(0)=0$
and $U(x)/|x|\to 0$ as $|x|\to \infty $.
With every $\pi=(B,C)\in \cA^x_a$ we associate the ``utility
process"
$$
J^{\pi}_t:=\int_0^{t\wedge \theta} e^{-\beta s} U(c_s)\,ds\,,\quad t\ge 0\,,
$$
where $\beta>0$. We consider the infinite horizon maximization
problem with  the {\em goal functional} $EJ^{\pi}_\infty$ and
define its {\it Bellman function} $W$ by \beq W(x):=\sup_{\pi \in
\cA^x_a}EJ^{\pi}_\infty\,,\quad x\in {\rm int}\,K\,. \eeq

Since $\cA^{x_1}_a\subseteq \cA^{x_2}_a$ when
$x_2-x_1\in K$, the function $W$ is {\bf increasing} with respect
to the partial ordering $\ge_K$ generated by the cone $K$.

If $\pi_i$, $i=1,2$, are admissible strategies for the initial
points $x_i$, then the strategy $\lambda \pi _1+(1-\lambda)\pi_2$
is an admissible strategy for the initial point $\lambda
x_1+(1-\lambda)x_2$, $\lambda \in [0,1]$, laying on the interval connecting $x_1$ and $x_2$. In the case where the relative price process $Y$
is continuous, the
corresponding ruin time for the process
\beq
\label{convprocess}
V^{\lambda x _1+(1-\lambda)x_2,\lambda \pi _1+(1-\lambda)\pi_2}=
\lambda V^{ x _1,\pi _1}+(1-\lambda)V^{x_2,\pi_2}
\eeq
 dominates the maximum of the ruin times
for  processes $V^{x_i,\pi_i}$.
The concavity of $u$ implies that
\beq
\label{convfunc}
J^{\lambda \pi _1+(1-\lambda)\pi_2}_t\ge \lambda J^{\pi _1}_t+(1-\lambda)J^{\pi_2}_t.
\eeq
and, hence, the function $W$ is
concave on ${\rm int}\,K$.

Unfortunately, in our main case of interest, where $Y$ has jumps,
the ruin times are not  related in such a simple way since the short positions are allowed.  It is easy to give examples of trajectories such that  $\theta^{ x _1,\pi _1} =\theta^{\lambda x _1+(1-\lambda)x_2,\lambda \pi _1+(1-\lambda)\pi_2}<\infty$ while $\theta^{ x _2,\pi _2}=\infty$ and the relations (\ref{convprocess}) and (\ref{convfunc}) do not hold. Therefore, we  cannot guarantee, by
the above argument, that the Bellman function is concave. Of course, these considerations show only that the concavity of $W$ cannot be obtained in a straightforward way as for a model based on continuous price process but  it is not excluded. 

The concavity of the Bellman function $W$ is not a property just
interesting {\it per se}. The classical definition of viscosity
solution, as was given by the famous ``User's guide" \cite{guide},
requires the continuity of $W$.  On the other hand, a concave function  is continuous in the interior of its domain (and even locally Lipschitz), see, e.g., \cite{Aubin}. Of course, the model must contains a provision which ensures that $W$
is finite. But the latter property, in the case of continuous price processes  implies, that $W$ is continuous on ${\rm int}\,K$.
In the case of processes with jumps one needs to analyze the continuity of $W$ using other arguments.

\smallskip
In the next section we show that the finiteness of $W$ still guarantees
its continuity in the interior of $K$.
 We do this using the following assertion.

\begin{lemm} Suppose that $W$ is a finite function. Let $x\in {\rm int}\, K$.
Then the function
$\lambda\mapsto W(\lambda x)$ is right-continuous on ${\bf R}_+$.
\end{lemm}
{\sl Proof.}
Let $\lambda>0$. Then
$\lambda\pi\in \cA_a^{\lambda x}$ if and only if  $\pi\in \cA_a^x$. For a concave function $U$ with $U(0)=0$ we have, for any $\e>0$ the inequality $U(c)\ge (1+\e)^{-1}U((1+\e)c)$.
Hence, for an arbitrary  strategy $\pi\in \cA_a^{x}$ we have, for $\theta=\theta^{x,\pi}=\theta^{(1+\e)x,(1+\e)\pi}$,  that
\bean
J_\infty^{(1+\e)\pi}-J_\infty^{\pi}&= &E\int_0^{\theta} e^{-\beta t}\big(U((1+\e)c_t)-U(c_t)\big)dt\\
&\le& \e E\int_0^{\theta} e^{-\beta t}U(c_t))dt\le \e W(x).
\eean
It follows that $W((1+\e)x)\le (1+\e)W(x)$. Since $W(x)\le W((1+\e)x)$, we infer from here that $\lambda\mapsto W(\lambda x)$
is right-continuous at the point $\lambda=1$. Replacing $x$ by $\lambda x$ we obtain the claim. \fdem

\smallskip
If $U$ is a homogeneous function of order $\gamma$ with $\gamma\in ]0,1[$, i.e.
$U(\lambda x)=\lambda^\gamma U(x)$ for all $\lambda>0$, $x\in K$,
then $W(\lambda x)=\lambda^\gamma W(x)$. Thus, the function
$\lambda\mapsto W(\lambda x)$ is concave and, therefore, continuous if finite.


\smallskip
\noindent{\bf Remark 1.} In financial models usually $\cC={\bf
R}_+e_1$ and $\sigma^0=0$, i.e. the only first (non-risky) asset
is consumed. Correspondingly, $U(c)=u(e_1c)=u(c^1)$ where $u$ is a utility function of a scalar argument.
Our presentation is oriented to the
power utility function $u_\gamma(x)=x^\gamma/\gamma$ with  $\gamma\in ]0,1[$. The case of $\gamma\le 0$, where, by convention, $u_0(x)=\ln x$, is of interest but it is not covered by the present study.

\smallskip
\noindent{\bf Remark 2.} We consider here a model with mixed
``regular-singular" controls. In fact, the assumption that the
consumption process has an intensity $c=(c_t)$ and the agent's
utility depends on this intensity is not very satisfactory from
the economical point of view. One can consider models with an
intertemporal substitution and the consumption by ``gulps", i.e.
dealing with ``singular" controls of the class $\cA^x$ and the
goal functionals like
$$
J_t^\pi:=\int_0^t e^{-\beta s} U(\bar C_s)ds\,,
$$
where
$$
\bar C_s=\int_0^s K(s,r)dC_r
$$
with a suitable kernel $K(s,r)$ (the exponential kernel
$e^{-\gamma (s-r)}$ is the common choice).

\section{Continuity of the Bellman Function}
\begin{prop}
\label{contBellman}
Suppose that $W(x)<\infty$ for all $x\in {\rm int}\, K$. Then $W$ is
continuous on ${\rm int}\, K$.
\end{prop}
{\sl Proof.}
First, we show that the function $W$ is upper semicontinuous on ${\rm int}\, K$.
Suppose that this is not the case and there is a sequence
$x_n$ converging to some $x_0\in {\rm int}\, K$ such that $\limsup_n W(x_n)>W(x_0)$. Without loss of generality we way assume that the sequence $W(x_n)$ converges.
The points $\tilde x_k=(1+1/k)x_0$, $k\ge 1$, belong to the ray ${\bf R}_+x_0$
and converges to $x_0$. We  find a subsequence $x_{n_k}$ such that  $\tilde x_k\ge_K x_{n_k}$ for all  $k\ge 1$.
Indeed, since
$$\tilde x_k=(1+1/k)x_0\in  x_0+{\rm int}\, K,
 $$
there exists $\e_k>0$ such that
 $$\tilde x_k+\cO_{\e_k}(0)\in  x_0+{\rm int}\, K.
 $$
It follows that
$$
\tilde x_k+(x_n-x_0)+\cO_{\e_k}(0)\in  x_n+{\rm int}\, K
$$
and, therefore, $\tilde x_k\in x_n+{\rm int}\, K$ for all $n$ such that
$|x_n-x_0|<\e_k$. Any strictly increasing sequence of indices $n_k$ with
$|x_{n_k}-x_0|<\e_k$ gives us in a subsequence of points $x_{n_k}$
 having the needed property. The function $W$ is increasing with respect
to the partial ordering $\ge_K$. Thus,
$$\lim_k W(\tilde x_k)\ge \lim_k W(x_{n_k})>W(x_0).
$$

On the other hand, the  function $\lambda\mapsto W(\lambda x_0)$ is right-continuous at $\lambda=1$ and, hence, $\lim_k W(\tilde x_k)=W(x_0)$.
This  contradiction shows that $W$ is upper semicontinuous on  ${\rm int}\, K$.

\smallskip

Let us show now that $\liminf_n W(x_n)\ge W(x_0)$ as $x_n\to x_0$, i.e. $W$ is lower semicontinuous on ${\rm int}\, K$.

Fix $\e>0$. Due to the finiteness of the Bellman function there are a strategy $\pi$ and $T\in {\bf R}_+$ such that for $\theta=\theta^{x_0,\pi}$ we have the bound
$$
E\int_0^{T\wedge \theta}e^{-\beta s}U(c_s)ds \ge W(x_0)-\e.
$$
It remains to show that
\beq
\label{thetan}
\liminf_n\, \theta_n\wedge T \ge \theta\wedge T\qquad \hbox{a.s.},
\eeq
where we use the abbreviation $\theta_n:=\theta^{x_n,\pi}$.
Indeed, with this bound we get, using the Fatou lemma, that
\bean
\liminf_n W(x_n)&\ge& \liminf_n E\int_0^{\theta_n\wedge T}e^{-\beta s}U(c_s)ds
\ge E\liminf_n \int_0^{\theta_n\wedge T}e^{-\beta s}U(c_s)ds
\\ &\ge& E\int_0^{\theta\wedge T}e^{-\beta s}U(c_s)ds\ge W(x_0)-\e
\eean
and the claim follows since $\e$ is arbitrarily small.

To prove (\ref{thetan}), we observe that by virtue of (\ref{dynam2}) on the interval $[0,\theta_n\wedge\theta\wedge T]$ we have the  representation
$$
V^{x_n,\pi}_t-V^{x_0,\pi}_t=D_{x_n-x_0}S_t
$$
implying that
$$
\sup_{t\le \theta_n\wedge\theta\wedge T}|V^{x_n,\pi}_t-V^{x_0,\pi}_t|\le S_T^*|x_n-x_0|,
$$
where $S_T^*:=\sup_{t\le T}|S_t|$.
Fix arbitrary, ``small", $\delta>0$. For almost all $\omega$ the distance $\rho(\omega)$ of the trajectory $V^{x_0,\pi}(\omega)$ from the boundary $\partial K$
on the interval $[0,\theta(\omega)\wedge T-\delta]$ is strictly positive. The above bound shows that for sufficiently large $n$ the  trajectory $V^{x_n,\pi}(\omega)$ does not deviate from $V^{x_0,\pi} (\omega)$ more than on $\rho(\omega)/2$ on the interval $[0,\theta_n(\omega)\wedge\theta(\omega)\wedge T]$.
It follows that $\theta_n(\omega)\ge \theta(\omega)\wedge T-\delta$. Thus,
$$
\liminf_n\, \theta_n\wedge T \ge \theta\wedge T-\delta\qquad \hbox{a.s.}
$$
and (\ref{thetan}) holds. \fdem

\section{The Hamilton--Jacobi--Bellman Equation}

 Let $G:=(-K)\cap \partial  \cO_1(0)$  where
 $\cO_r(y):=\{x\in {\bf R}^d:\
 |x-y|< r\}$.  The set $G$ is  compact   and $-K={\rm
cone}\, G$. We denote by $\Sigma_G$  the {\em support function} of
$G$, given by the relation $\Sigma_G(p)=\sup_{x\in G}px$.
The convex function $U^*(.)$ is the Fenchel dual of the convex function $-U(-.)$ whose domain is $-\cC$, i.e.
$$
U^*(p)=\sup_{x\in \cC}(U(x)-px).
$$

We denote by $C_1(K)$ the subspace of the space of continuous functions  $f$ on $K$ such that $\sup_{x\in K}|f(x)|/(1+|x|)<\infty$. 
In other words, $C_1(K)$ is the space of continuous functions on $K$ of sublinear growth.  The notation $f\in C^2(x)$ means that $f$ is smooth in some neighborhood of $x$. 
\smallskip

Let $f\in C_1(K)\cap C^2({\rm int}\, K)$.
Using the abbreviation 
$$
I(z,x):=I_{\{z:\ x+D_xz \in {\rm int}\, K\}}=I_
{{\rm int}\, K}( x+D_xz)
$$
we introduce  the function
$$
\cI(f,x): =\int \big[(f(x+D_xz)I(z,x)-f(x))-D_xz f'(x)\big]\Pi(dz), \quad x\in {\rm int}\, K.
$$
It is well-defined and continuous in $x$. Indeed,  fix $x_0\in {\rm int}\, K$. Let $\e\in ]0,1]$ be such that the ball $\cO_{4\e}(x_0)\subset K$.  With this choice
$x+D_xz\in \cO_{2\e}(x_0)$ when $x\in \cO_{\e}(x_0)$ and $|z|\le \delta:=\e/(1+|x_0|)$.
Using the Taylor formula for such value of $z$ and
the sublinear growth of $f$ for $z$ with $|z|> \delta$  we obtain the following uniform bound for $x\in \cO_{\e}(x_0)$:
$$
|(f(x+D_xz)I(z,x)-f(x)-D_xzf'(x)|\le \kappa_1 |z|^2I_{\cO_{\delta}(x_0)}(z) +\kappa_2|z|I_{\cO^c_{\delta}(0)}(z).
$$
It implies the needed integrability and the continuity of the integral in $x$.

 We introduce a  function of five variables
by putting
$$
F(X,p,\cI(f,x),W,x):=\max \{ F_0(X,p,\cI(f,x),W,x)+U^*(p), \Sigma_G(p)\},
$$
where $X$ belongs to $\cS_d$, the set of $d\times d$  symmetric
matrices, $p,x\in {\bf R}^d$, $W\in {\bf R}$, $f\in C_1(K)\cap C^2(x)$, and the function
$F_0$ is given by
$$
F_0(X,p,\cI(f,x),W,x):=\frac 12 {\rm tr}\, A(x) X+ \mu (x)p+\cI(f,x)-\beta W
$$
where $A(x)$ is the matrix with  $A^{ij}(x):=a^{ij}x^ix^j$, and $\mu(x)$ is the vector with components  $\mu ^{i}(x):=\mu ^{i}x^i$, $1\le i,j\le d$.

In a more detailed form we have that
$$
F_0(X,p,\cI(f,x),W,x)=\frac 12 \sum_{i,j=1}^d a^{ij}x^ix^jX^{ij}+
\sum_{i=1}^d\mu^ix^ip^i+\cI(f,x)-\beta W.
$$

Note that $F_0$ is increasing in the argument $f$ in the same sense as $\cI$.

 If $\phi $ is a smooth function, we put
$$
\cL\phi(x):=F(\phi''(x),\phi'(x),\cI(\phi,x), \phi(x),x).
$$
In a similar way,   $\cL_0$ corresponds to the function $F_0$.

We show, under mild hypotheses, that $W$ is a  viscosity
solution of the Dirichlet problem for the HJB equation \bea
\label{vi} F(W''(x),W'(x),\cI(W,x),W(x),x)&=&0,
 \qquad x\in  {\rm int}\, K,\\
W(x)&=&0,
 \qquad x\in  \partial K,
 \label{icond}
\eea with the boundary condition understood in the usual classical
sense and establish a uniqueness result for this problem.

\section{Viscosity Solutions for  Integro-Differential Operators}

Since, in general, $W$ may have no  derivatives at some points
$x\in {\rm int} K$ (and this is, indeed, the case for the model
considered here), the notation (\ref{vi}) needs to be interpreted.
The idea of viscosity solutions is to substitute $W$ in $F$ by suitable test functions. Formal definitions (adapted to the case we
are interested in) are as follows.

A  function $v\in C(K)$ is called  {\it viscosity supersolution}
of (\ref{vi}) if  for every  $x\in {\rm int}\, K$ and  every   $f\in C_1(K)\cap C^2(x)$ such that
$v(x)=f(x)$ and $v \ge f$
 the inequality $\cL f(x)\le 0$ holds.

A  function $v\in C(K)$  is called {\it viscosity subsolution} of (\ref{vi})
 if  for every  $x\in {\rm int}\, K$ and every  $f\in C_1(K)\cap C^2(x)$ such that
$v(x)=f(x)$ and $v \le f$
 the inequality $\cL f(x)\ge 0$ holds.

A  function $v\in C(K)$  is a {\it viscosity solution} of
(\ref{vi}) if $v$ is simultaneously a viscosity super- and
subsolution.

At last, a  function $v\in C_1(K)\cap C^2({\rm int}\, K)$  is called  {\it classical
supersolution} of (\ref{vi}) if  $\cL
v \le 0$ on ${\rm int}\, K$. We add the adjective {\it strict}
when $\cL v <0$ on the set ${\rm int}\, K$.


\smallskip
For the sake of simplicity and having in mind the specific case we
shall work on, we incorporated in the definitions the requirement
that the viscosity super- and subsolutions  are continuous on $K$
including the boundary. For other cases this might be too
restrictive and more general and flexible formulations can be
used.

\begin{lemm}\label{classicsol}
Suppose that the function $v$ is a viscosity solution of
(\ref{vi}). If $v$ is twice differentiable at  $x_0\in {\rm int}\, K$, then it
satisfies (\ref{vi}) at this point in the classical sense.
\end{lemm}
 {\sl Proof.} One needs to be more precise with definitions
 since it is not assumed
 that $v'$ is  defined at every point of a neighborhood of $x_0$.
 ``Twice differentiable" means here that  the Taylor
 formula at $x_0$ holds:
 $$
v(x)= P_2(x-x_0) + (x-x_0)^2h(|x-x_0|)
 $$
 where
 $$
 P_2(x-x_0):=v(x_0)+\langle v'(x_0),x-x_0 \rangle +\frac 12 \langle
v''(x_0)(x-x_0),x-x_0\rangle
 $$
 and  $h(r)\to 0$ as $r\downarrow 0$.
We introduce the notation
$\Gamma_r:=\{z\in {\bf R}^d:\ |{\rm diag}\,x_0 z|\le r\}$, $r> 0$. Note that $\cO_{r/|x_0|}\subset \Gamma_r$. Hence, $\Pi(\Gamma^c_r)<\infty$.

Let $\e\in]0,1]$. We choose a number
 $\delta_0 \in ]0,1[$ such that $x_0+\cO_{\delta_0}(0)\subset {\rm int}\, K$ and $|h(s)|\le\e $ for $s\le \delta_0$. Put   $\delta:=\delta_0/(1+|x_0|)$. Take $\Delta\in]\delta,\delta_0[$ sufficiently close to $\delta$ to insure that
 $x_0+\cO_\Delta(0)\subset {\rm int}\, K$  and  $\Pi(\Gamma_\Delta\setminus \Gamma_\delta)\le\e$.

We define the function $f_\e\in C_1(K)\cap C^2(x_0) $ by the formula
 $$
f_\e(x)=\left \{
\begin{array}{ll}
P_2(x-x_0) + \e(x-x_0)^2,&\quad  x\in x_0+\cO_\delta(0),\\
g(x)\vee v(x),&\quad x\in x_0+\cO_\Delta(0)\setminus \cO_\delta(0),\\
v(x),& \quad x\in x_0+\cO^c_\Delta(0),
\end{array}\right.
 $$
where
$$
g(x):=P_2\left (\delta \frac{x-x_0}{|x-x_0|}\right) + \e\delta + \frac {\delta - |x-x_0|}{\Delta-|x-x_0|}.
$$
Clearly,  $f_\e(x_0)=v(x_0)$ and  $f_\e\ge v$.
Since $v$ is a viscosity subsolution,  $\cL f_\e(x_0)\ge 0$.
Note that
$$
|\cL f_\e(x_0)-\cL v(x_0)|\le \e\sum_{i=1}^na^{ii}(x^i_0)^2+\cI(f_\e-v,x_0),
$$
with
$$
\cI(f_\e-v,x_0)=\int (f_\e-v)(x_0+D_{x_0} z)I_{\{x_0+D_{x_0} z\in {\rm int}\, K\}} \Pi(dz).
$$
Let us check that $\cI(f_\e-v,x_0)\le \kappa \e$. Indeed,
\bean
(f_\e-v)(x_0+D_{x_0} z)&\le& \e(D_{x_0} z)^2I_{\Gamma_\delta}+M I_{\Gamma_{\Delta}\setminus \Gamma_\delta}\\
&\le&
\e\min\{|x_0|^2|z|^2,\delta^2\}+M I_{\Gamma_{\Delta}\setminus \Gamma_\delta}.
  \eean
where $M=1+\sup_{y\in\cO_1(0)}|P_2(y)|$. It follows that
$$
\cI((f_\e-v) I_{\cO^c_\Delta(0)},x_0)\le \e(1+|x_0|)^2\int|z|^2\wedge 1\Pi(dz)+M\e.
$$
 Letting $\e$ tend to zero, we
obtain that $\cL v(x_0)\ge 0$. Arguing similarly with $\e<0$, we get the opposite
inequality. \fdem

\section{Jets}
Let $f$ and $g$ be functions defined in a neighborhood of zero. We
shall write $f(.)\lessapprox g(.)$ if $f(h)\le g(h)+o(|h|^2)$ as
$|h|\to 0$. The notations $f(.)\gtrapprox g(.)$ and $f(.)\approx
g(.)$ have the obvious meaning.

For $p\in {\bf R}^d$ and $X\in \cS_d$ we consider the quadratic
function
$$Q_{p,X}(z):=pz+(1/2)\langle Xz,z\rangle\,,\quad z\in {\bf R}^d\,,
$$
and define  the {\em super}- and {\em subjets} of a function $v$
at the point $x$: \bean
J^+v(x)&:=&\{(p,X):\ v(x+.)\lessapprox v(x)+Q_{p,X}(.)\},\\
J^-v(x)&:=&\{(p,X):\ v(x+.)\gtrapprox v(x)+Q_{p,X}(.) \}. \eean

In other words, $J^+v(x)$ (resp. $J^-v(x)$) is the family of
coefficients of quadratic functions $v(x)+Q_{p,X}(y-.)$ dominating  the
function $v(.)$ (resp., dominated by this function) in a neighborhood of the point $x$
with  precision up to the second order included and coinciding
with $v(.)$ at this point.

  In the classical theory developed for differential equations the notion of viscosity solutions admits an equivalent formulation in terms of super- and subjets. Since the latter are ``local" concepts, such a characterization is not possible for integro-differential operators. Nevertheless, one can construct from semijets test functions with useful properties.

The following lemma claims for $v\in C_1(K)$   with any element  $(p,X)\in J^+v(x)$  $x\in {\rm int}\,K$ one can relate a test function dominating $v$, arbitrary close to $v$ in the uniform metric,  touching $v$ at the point $x$, smooth at a neighborhood of  $x$ and having at this point the first and the second derivatives coinciding with $p$ and $X$.     


\begin{lemm}
 \label{jet}Let $v\in C_1(K)$ and let $\alpha>0$. Let $x\in {\rm int}\,K$ and let  $(p,X)\in J^+v(x)$.  
Then there exist a number $a_0\in ]0,1[$ and a $C^2$-function $r: {\bf R}^d\to {\bf R}$ with compact support such that 
\beq
\label{jet}
 \lim_{|h|\to 0}|u|^{-2}r(h)=0,  
\eeq
 and the function $f_0:K\to {\bf R}$ 
given by the formula 
\beq
\label{strf}
f_0(x+h):=\big[\big(v(x)+Q_{p,X}(h)+r(h)\big)\vee v(x+h)\big]\wedge (v(x+h)+\alpha),  \quad x+h\in K,
\eeq
has the following properties:   
$$
f_0(x+h)=v(x)+Q_{p,X}(h)+r(h), \qquad h\in \cO_{a_0}(0), 
 $$ 
$v\le f_0\le v+\alpha $ on $K$,  $f_0(x)=v(x)$, $f'_0(x)=p$, $f''_0(x)=X$. 

  In particular, if $v$ is a subsolution of the HJB equation, then 
$\cL f\le 0$ on ${\rm int}\, K$. 
\end{lemm}
\noindent
{\sl Proof.} Take $a_0\in ]0,1[$ such that the ball $\cO_{2a_0}(x)=\{y\in {\bf R}^d:\ |y-x|\le 2a\}$ lays in the interior of $K$.
By definition,
$$
v(x+h)-v(x)-Q_{p,X}(h)\le |h|^2 \varphi (|h|),
$$
where $\varphi(u)\to 0$ as $u\downarrow 0$. We consider on $]0,a_0[$
the  function
$$
\delta(u):=\sup_{\{h:\ |h|\le u\}}\frac 1{|h|^{2}}
(v(x+h)-v(x)-Q_{p,X}(h))^+\le \sup_{\{y:\ 0\le y \le u\}}\varphi^+
(y). 
$$
Obviously, $\delta $ is continuous, increasing and $\delta(u)\to
0$ as $u\downarrow 0$. We extend $\delta $ to a continuous function on ${\bf R}_+$
with $\delta(u)=0$ for $u\ge 1$. 

The function
$$
\Delta(u):=\frac 2{3}\int _u^{2u}\int _\eta^{2\eta}\delta(\xi)d\xi
d\eta
$$
vanishes at zero with its two right derivatives; $u^2\delta(u)\le
\Delta(u)\le u^2\delta(4u)$. It follows that the function
$r: h\mapsto \Delta(|h|)$ has a compact support, belongs to $C^2(\cO_{a_0}(0))$, its Hessian
vanishes at zero, and
$$
v(x+h)-v(x)-Q_{p,X}(h)\le |h|^2\delta(|h|)\le \Delta(|h|)=r(h), \qquad h\in \cO_{a_0}(0)
$$

Thus, the function $y\mapsto v(x)+Q_{p,X}(y-x)+r (y-x)$ dominates 
$v$ on the ball $\cO_{a_0}(x)$. Without loss of generality, diminishing $a_0$ if necessary, we may assume that it is dominated by    
$v + \alpha$ on this ball. Now the assertion of the lemma is obvious. 
\fdem

\smallskip
The corresponding assertion for $J^-v(x)$  also holds --- with obvious changes in the formulation.  

\smallskip
For the proof of the uniqueness theorem we need specific families of test functions coinciding with sub- and supersolutions  outside a neighborhood of $x$. 
To this end we introduce the following definitions. 

Let $0<a<a'$. We say that a continuous mapping $\xi_{a,a'}: {\bf R}^d\to [0,1]$ is an {\it $(a,a')$-cutoff function} if  $\xi_{a,a'}=1$
on $\cO_a(0)$ and $\xi_{a,a'}=0$
outside  $\cO_{a'}(0)$. If $L$ is a linear subspace of ${\bf R}^d$ we define the {\it cylindrical $(a,a')$-cutoff function}  $\xi^L_{a,a'}$ by putting $\xi^L_{a,a'}(x)=\xi_{a,a'}(P_Lx)$
where $P_L$ is a projection of $x$ onto $L$. 

It is clear that in the notation of the above lemma for any $a'\in ]0,a_0]$ the functions $f:K\to {\bf R}$ given by the formulae 
\bea
\label{strf}
f(x+h)&:=&\big(v(x)+Q_{p,X}(h)+r(h)\big)\xi_{a,a'}(h)+v(x+h)\big(1-\xi_{a,a'}(h)\big), \quad x+h\in K,\\
\label{cstrf}
f(x+h)&:=&f_0(x+h)\xi^L_{a,a'}(h)+v(x+h)\big(1-\xi^K_{a,a'}(h)\big), \quad x+h\in K,
\eea
will satisfy all the properties claimed for $f_0$.

The following lemma will be used in the case when $D=D_x={\rm diag}\, x$, $\tilde D=D_y={\rm diag}\, y$, and $x$, $y$ has no zero components. 

\begin{lemm} \label{Mollif=} 
Let $D,\,\tilde D$ be two invertible linear operators on ${\bf R}^d$. Let $\xi_{a,a'}$ be a $(a,a')$-cutoff function. Then there is a  $(\tilde a,\tilde a')$-cutoff function $\tilde \xi_{\tilde a,\tilde a'}$ with arbitrary small  $\tilde a\le a||D\tilde D^{-1}||^{-1}$ and arbitrary  $\tilde a'\ge a'||\tilde DD^{-1}||$ such that 
$$
\tilde \xi_{\tilde a,\tilde a'}(\tilde Dz)=\xi_{a,a'}(D z)\quad \forall\,z\in \mathbf{R}^d.
$$
\end{lemm}
{\sl Proof.} Put
$$
\tilde \xi_{\tilde a,\tilde a'}(u):=\xi_{a,a'}(D\tilde D^{-1}u), \quad u\in \mathbf{R}^d.   
$$
Then $\xi_{\tilde a,\tilde a'}(u)=1$ if $|u|\le a||\tilde D D^{-1}||^{-1}$, and $\tilde \xi_{\tilde a,\tilde a'}(u)=0$ if $|D\tilde D^{-1}u|\ge a'$. The last inequality holds when $|u|\ge a'||(D\tilde D^{-1})^{-1}||=a'||\tilde DD^{-1}||$. \fdem 

\begin{remark} The assertion of the lemma remains true for cylindrical cut-off functions   in the situation where $D$ and $\tilde D$ are two symmetric operators with 
the common image space $L={\rm Im}\, D={\rm Im}\, \tilde D$. The norm in formulation is the norm in $L$ of their restrictions. 
\end{remark}

\section{Supersolutions and Properties of the Bellman Function}
\subsection{When  is the Bellman Function $W$ Finite on $K$?}

First, we present sufficient conditions ensuring that the Bellman
function   $W$ of the considered maximization problem is finite.

Functions we are interested  in are defined in the solvency cone $K$ while the process $V$  may jump out of the latter.
 In order to be able to apply later the It\^o formula we stop $V=V^{x,\pi}$ at the moment immediately
preceding the ruin and define the process
  $$\tilde V=V^{\theta-}=VI_{[0,\theta[}+V_{\theta-} I_{[\theta,\infty[},$$
where $\theta$ is the exit
time  of  $V$ from the interior of the solvency cone $K$. This process coincides with $V$ on
$[0,\theta[$ but, in contrast to the latter, either always remains
in $K$ (due to the stopping at $\theta$ if
$V_{\theta-}\in {\rm int}\, K$) or exits to the boundary in a
continuous way and stays on it at the exit point.

Since $\tilde V_t=V_{t\wedge \theta}-\Delta V_\theta I_{[\theta,\infty[}(t)$, we obtain from the equation (\ref{dynamvec}) the representation
\bean\label{tildeV}
\tilde V_t&=&x+\int_0^{t\wedge \theta}D_{\tilde V_{s-}} (\mu ds+
\Xi dw_s)
+\int_0^{t\wedge \theta} \int D_{\tilde V_{s-}}z(p(dz,ds)-q(dz,ds))
\\
&&-\Delta V_\theta I_{[\theta,\infty[}(t)
+B_t-C_t.
\eean

Let $\Phi$ be the set of  continuous functions  $f:K\to {\bf R}_+$
increasing with respect to the partial ordering $\ge_K$ and such
that for every   $x\in {\rm int}\, K$ and  $\pi\in \cA^x_a$ the
positive process $X^f=X^{f,x,\pi}$ given by the formula
 \beq
\label{Xf} X^f_t:= e^{-\beta t}f(\tilde V_t)I_
{[0,\theta[}({t}) + J^{\pi}_t
 \eeq
 is a supermartingale.

 The set $\Phi$ of $f$ with
this  property is convex and stable under the operation $\wedge$
(recall that the minimum of two supermartingales is a
supermartingale). Any continuous function which is a monotone
limit (increasing or decreasing) of functions from $\Phi$ also
belongs to $\Phi$.

The interest to the processes $X^f$ with $f\in \Phi$ is explained by the following:  

\begin{lemm}
\label{Phi} $(a)$ If $f\in \Phi$, then $W\le f$.

$(b)$ Let  $y\in \partial K$. Suppose that for every $\e>0$  there exists  $f_\e\in \Phi$ such
that $f_\e(y)\le \e$.  Then $W$ is continuous at $y$ and $W(y)=0$. 
\end{lemm}
{\sl Proof.} $(a)$ On the boundary $\partial K$ the inequality is trivial. 
Using the positivity of  $f$, the
supermartingale property of $X^f$, and, finally, the monotonicity
of $f$ we get, for $\in {\rm int}\, K$, the following chain of inequalities leading to the
required property:
$$
EJ^{\pi}_t\le EX^f_t\le f(\tilde V_0)=f(V_0)\le f(V_{0-})= f(x).
$$

$(b)$  The
continuity of the function $W$ at the point $y\in
\partial K$ follows from the inequalities
 $0\le W\le f_\e$.  
\fdem

\smallskip

\noindent
{\bf Remark.}  Recall that 
Proposition \ref{contBellman} assets that the function $W$, if finite,  is continuous 
on  the interior of $K$. Thus, the above lemma implies that $W$ is  continuous  on
 ${\rm int}\, K$ if  $\Phi$ is not empty.  If $\Phi$ is reach enough to apply $(b)$ at every point of the boundary, then $W$ is continuous on $K$ and vanishes on the boundary. 
 \smallskip

\begin{lemm}
\label{lemmsuper} Let $f:K\to {\bf R}_+$ be a function in
$C_1(K)\cap C^2({\rm int}\, K)$. If $f$ is a classical supersolution
of\, (\ref{vi}), then $f\in \Phi$, i.e. $f$ is 
increasing with respect to the partial ordering $\ge_K$ and  $X^f$ is a
supermartingale.
\end{lemm}
{\sl Proof.} First, notice that a classical supersolution is
 increasing with respect to the
partial ordering $\ge_K$. Indeed, by the finite increments formula
we have that for any $x,h\in {\rm int}\, K$
$$ f(x+h)-f(x)=f'(x+\vartheta h)h
$$
for some $\vartheta\in [0,1]$. The right-hand side is greater or
equal to zero because for the supersolution $f$ we have the
inequality $\Sigma_G(f'(y))\le 0$ whatever is $y\in {\rm int}\,
K$, or, equivalently, $f'(y)h\ge 0$ for every $h\in K$, just by
the definition of the support function $\Sigma_G$ and the choice
of $G$  as a generator of the cone $-K$. By continuity,
$f(x+h)-f(x)\ge 0$ for every $x,h\in K$.

Let  $\theta_n:=\inf\,\{t\colon\  {\rm dist}(\tilde V_t, \partial K)\le 1/n \} $. The stopped processes 
$\tilde V^{\theta_n}$ evolves in ${\rm int}\, K$. Thus, we can apply the 
 ``standard" It\^o formula to $e^{-\beta t}f(\tilde V_t)$ and obtain, for $t\le \theta$, 
that
\bean
e^{-\beta t}f(\tilde V_t)&=&f(x)+\int_0^t e^{-\beta s}f'(\tilde V_{s-})d\tilde V_{s}  -\beta
\int_0^t e^{-\beta s}f(\tilde V_{s-})ds \\
&&
+ \frac 12 \int_0^t e^{-\beta s} {\rm tr}\, A(\tilde V_{s-})f''(\tilde V_{s-})ds \\
&&+\sum_{s\le t}e^{-\beta s} [f(\tilde V_{s-}+\Delta
\tilde V_s)-f(\tilde V_{s-})-f'(\tilde V_{s-})\Delta\tilde V_{s}].
\eean

Taking into account that the processes $Y$ and $B$ do not jump simultaneously and the 
ruin cannot happen due to a jump of $B$ we get that   
{\small
\bean
&&\sum_{s\le t}e^{-\beta s} [f(\tilde V_{s-}+\Delta
\tilde V_s)-f(\tilde V_{s-})-f'(\tilde V_{s-})\Delta\tilde V_{s}]-e^{-\beta\theta}f'(V_{\theta-})\Delta V_\theta I_{\{\theta\}}(t)\\
&&-e^{-\beta\theta}f(V_{\theta-})I_{\{\theta\}}(t)
\\
&=&\sum_{s\le t}e^{-\beta s} [f(V_{s-}+\Delta
V_s)I_
{{\rm int}\, K}( V_{s-}+\Delta V_{s})-f( V_{s-})-f'(V_{s-})\Delta  V_{s}]\\
&=&\int_0^t\int e^{-\beta s} [f(V_{s-}+D_{V_{s-}}z)I(V_{s-},z)-f( V_{s-})-f'(V_{s-})D_{
V_{s-}}z)]I_{\{\Delta B_s=0\}}p(ds,dz)\\
&&+\sum_{s\le t} e^{-\beta s} [f( V_{s-}+\Delta
B_s)-f(V_{s-})-f'(V_{s-})\Delta B_s]\\
&=&
\int_0^t\int e^{-\beta s} [...]I_{\{\Delta B_s=0\}}(p(ds,dz)-\Pi(dz)ds)+\int_0^t\int e^{-\beta s} [...]\Pi(dz)ds\\
&&+\sum_{s\le t} e^{-\beta s} [f( V_{s-}+\Delta
B_s)-f(V_{s-})-f'(V_{s-})\Delta B_s]
\eean
}
where we replace in the integrals  by dots the lengthy expression $$
f(\tilde V_{s-}+D_{\tilde V_{s-}}z)I(V_{s-},z)-f(\tilde V_{s-})-f'( \tilde V_{s-})D_{\tilde V_{s-}}z.
$$
Noting that 
$$
X^f_t=e^{-\beta t}f(\tilde V_t)-e^{-\beta \theta}f(V_{\theta-})I_{\{\theta\}}(t)+J^{\pi}_t
$$
and using the equation (\ref{dynamvec}) and the above formulae we obtain, after regrouping terms, the following representation 
\beq \label{Ito}
X^f_t=f(x)+\int_0^{t\wedge \theta} e^{-\beta s} [\cL_0 f(\tilde V_s) -c_s
f'(\tilde V_s)+U(c_s)]ds +R_t+ m_t,
\eeq
where
 \beq \label{R} R_t:=\int_0^{t\wedge \theta} e^{-\beta s}f'(
V_{s-})dB^c_s+\sum_{s\le t} e^{-\beta s} [f(\tilde  V_{s-}+\Delta
B_s)-f(\tilde  V_{s-})]
\eeq
and $m$ is the local martingale
{\small
\bea
\nonumber
m_t&=&\int_0^{t\wedge \theta}e^{-\beta s}
f'(\tilde V_{s-})D_{\tilde V_{s-}}
\Xi dw_s\\ 
&&+\int_0^{t\wedge\theta}\int e^{-\beta s} [f(\tilde V_{s-}+D_{\tilde V_{s-}}z)I(\tilde V_{s-},z)-f(\tilde V_{s-})](p(dz,ds)-\Pi(dz)ds). \label{lm}
\eea
}
By definition of a supersolution,
 for any $x\in {\rm int}\, K$,
$$
\cL_0 f(x)\le - U^*(f'(x))\le cf'(x)-U(c)\quad \forall\,c\in \cC.
$$
Thus, the integral in (\ref{Ito}) is a decreasing process. The
process $R$ is also decreasing.  Indeed, the terms of the sum in
(\ref{R}) are less or equal to zero in virtue of the monotonicity of $f$ and
$$
f'(V_{s-})dB^c_s=I_{\{\Delta B_s=0\}}f'( V_{s-})\dot
B_s d||B||_s
$$
where $f'(V_{s-})\dot B_s\le 0$ since $\dot B$ takes values
in $-K$.
Let $\sigma_n$ be a localizing sequence for $m$.
Taking into account that $ X^f\ge 0$, we obtain from
(\ref{Ito}) that for each $n$ the negative decreasing process $R_{t\wedge
\sigma_n}$ dominates an integrable process and so it is integrable.
The same conclusion holds for the stopped integral. Being a sum of
an integrable decreasing process and a martingale, the process
$ X^f_{t\wedge \sigma_n}$ is a positive supermartingale and, hence,
by the Fatou lemma, $ X^f$ is a supermartingale as well.
\fdem


\smallskip

Lemma \ref{lemmsuper} implies that the existence of a smooth
positive supersolution $f$ of (\ref{vi}) ensures the
finiteness of $W$ on $K$.
We discuss a method how to construct  supersolutions in Section 12.

\smallskip

\noindent {\bf Remark.}
 Let $\bar \cO$ be the closure of an open subset $\cO$ of $K$ and let
$f:\bar \cO\to {\bf R}_+$ be a classical supersolution in $\bar \cO$ increasing with respect to  the partial ordering $\ge_K$. Let
$x\in \cO$ and let $\tau$ be the exit time of the
process $V^{x,\pi}$ from $\bar \cO$. The above arguments imply that the
 process $X^f_{t\wedge \tau}$ is a supermartingale and, therefore,
\beq
\label{localss} E[e^{-\beta (t\wedge \tau)}f(\tilde V_{t\wedge \tau})I_
{[0,\theta[}({t\wedge \tau}) +
J^{\pi}_{t\wedge \tau}]\le f(x).
\eeq

\subsection{Strict Local Supersolutions}

For the strict supersolution we can get a more precise result  which
will play the crucial role in deducing from the Dynamic
Programming Principle the property of $W$ to be a subsolution of the
HJB equation.

Fix $x\in {\rm int}\,K$  and a ball $\bar \cO_r(x)\subseteq {\rm int}\,K$ such that the larger ball $\bar\cO_{2r}(x)\subseteq {\rm int}\,K$.  We define
$\tau^\pi=\tau^\pi_r$ as the exit time of $V^{\pi,x}$ from $ \cO_r(x)$, i.e.
$$
\tau^\pi:=\inf\{t\ge 0: \ |V^{\pi,x}_t-x|\ge r\}.
$$

\begin{lemm}
\label{svi} Let $f\in C_1(K)\cap C^2(\cO_{2r}(x))$ be such that $\cL f\le
-\e<0$ on $\bar \cO_r(x)$. Then there exist a constant $\eta=\eta_\e>0$
and an interval $]0,t_0]$ such that
$$ \sup_{\pi\in
\cA^{x}_a} E X^{f,x,\pi}_{t\wedge \tau^\pi}\le f(x)- \eta t \qquad
\forall\,t\in ]0,t_0].
$$
\end{lemm}
{\sl Proof.} We fix a strategy $\pi$ and omit its symbol in the
notations below. In what follows, only the behavior of the
processes on $[0,\tau]$ does matter. Note that $|V_\tau-x|\ge r$ on
the set $\{\tau<\infty\}$ and $\tau\le \theta$. As in the proof of Lemma
\ref{lemmsuper}, we apply the It\^o formula and obtain, with the same notations, 
(\ref{R}) and (\ref{lm}), 
the representation
\bean
X^f_{t\wedge \tau}&:=&e^{-\beta (t\wedge \tau)}f(\tilde V_{t\wedge \tau})I_
{[0,\theta[}({t\wedge \tau})+J^{\pi}_{t\wedge \tau}\\
&=&f(x)+
\int_0^{t\wedge \tau} e^{-\beta s} (\cL_0 f+U^*)(\tilde V_s) ds \\
&& -
\int_0^{t\wedge \tau} e^{-\beta s} [U^*(\tilde V_s) +c_s
f'(\tilde V_s)-U(c_s)]ds
+R_{t\wedge \tau}+ m_{t\wedge \tau}.
\eean
Due to the monotonicity of $f$ we may assume
without loss of generality that on the interval $[0,\tau]$ the increment $\Delta B_t$ does not exceed  the distance of $V_{s-}$ to the boundary of $ \cO_r(x)$. In other words, if the exit from the ball is due to an  action (and not because of a jump of the price process),  we can replace this action by a less expensive one,  with the jump of the process $\tilde V$ in the same direction but a smaller one, ending on the boundary of the ball. So, $|\Delta B_t|\le 2r$ for $t\le \tau$.  

By assumption, for $y\in \bar \cO_r(x)$ we have the bounds  $ (\cL_0 f+U^*)(y)\le
-\e$ (implying that the first integral in the right-hand side above is dominated by 
$- \e \; (t\wedge \tau)$) and $\Sigma_G(f'(y))\le -\e$. The latter inequality means that
$kf'(y)\le -\e|k|$ for every $k\in -K$ (therefore, we have the inclusion $f'(\bar \cO_r(x))\subset
{\rm int}\,K^*$). In particular, for  $s\in [0,\tau]$
$$
f'(V_{s-})\dot B_s\le  -\e|\dot B_s|, \qquad [f(\tilde  V_{s-}+\Delta
B_s)-f(\tilde  V_{s-})]\le  -\e|\Delta B_s|.
$$
Since $|\tilde V_{s-}-x|\le r$ for $s\in[0,\tau]$, we obtain,
using the finite increment formula and the linear growth of $f$, the bounds
\bean
[f(\tilde V_{s-}+D_{\tilde V_{s-}}z)-f(\tilde V_{s-})]^2I(\tilde V_{s-},z)I_{\{|z|\le 1/2\}}&\le &\kappa |z|^2I_{\{|z|\le 1/2\}},
\eean
\bean
[f(\tilde V_{s-}+D_{\tilde V_{s-}}z)-f(\tilde V_{s-})]I(\tilde V_{s-},z)I_{\{|z|> 1/2\}}&\le& \kappa (1+|z|)I_{\{|z|> 1/2\}}, 
\eean
and, as  $I(\tilde V_{s-},z)=1$ when $|z|< r/(|x|+r)$,
$$
f(\tilde V_{s-})(1-I(\tilde V_{s-},z))\le \kappa I_{\{|z|\ge r/2\}} 
$$
for some constant $\kappa$ independent on the strategy. Thus, the integrand in the stochastic integral with respect to the centered Poisson measure in (\ref{lm}) for $t\le\tau$ is bounded by the function $|z|^2\wedge|z|$ multiplied by a constant while the integrand in the integral 
with respect to the Wiener process is bounded.  
It follows
that the local martingale $(m_{t\wedge\tau})_{t\ge 0}$  is a martingale and $Em_{t\wedge\tau}=0 $.

The above observations imply the inequality
$$
EX_{t\wedge \tau}^{f,x}\le f(x) - e^{-\beta t}E N_{t},
$$
where
$$
N_{t}:= \e \; (t\wedge \tau)+ \int _0^{t\wedge \tau}H(c_s,f'(V_s))ds
+\e \int _0^{t\wedge \tau}|\dot B_s|d||B||_s
$$
with $H(c,p):=U^*(p)+pc-U(c)\ge 0$. It remains to verify that
$EN_t$ dominates, on a certain interval $]0,t_0]$,
 a strictly increasing linear function which is independent of
 $\pi$.

The process $N_t$ looks a bit complicated but we can  replace it by another one of a simpler structure. To this end, note that there is a constant $\kappa$ (``large", for convenience,
  $\kappa \ge 1$)  such that
$$
\inf_{p\in  f'(\bar \cO_r(x))} H(c,p)\ge \frac{\e}{2} |c|, \qquad
\forall\,  c\in \cC , \quad|c|\ge \kappa.
$$
Indeed, being the image of a closed ball under continuous mapping, the set
$f'(\bar \cO_r(x))$ is a compact in ${\rm int}\,K^*$. The lower
bound of the continuous function $U^*$ on $f'(\bar \cO_r(x))$ is finite.
 For any $p$ from $f'(\bar \cO_r(x))$ and $c\in \cC\subseteq K $ we have
the inequality $pc/|c|\ge \e$.  At last,
 $U(c)/|c|\to 0$ as $c\to \infty$. Combining these facts we infer
 the claimed inequality.
Thus, for the first integral in the definition of $N_t$ we have the bound
$$
\int _0^{t\wedge \tau}H(c_s,f'(V_s))ds\ge \frac{\e}{2}
\int_0^{t\wedge \tau}I_{\{|c_s|\ge \kappa\}} |c_s|ds.
$$
The second integral in the definition dominates  $ \kappa_1
||B||_{t\wedge \tau}$ for some $\kappa_1>0$. To see this, let us consider the absolute norm $|.|_1$ in ${\bf R}^d$. In contrast with the total variation $||B||$ which is calculated with respect to the Euclidean norm $|.|$, the total variation of $B$ with respect to the absolute norm admits a simpler expression  $\sum_i{\rm
Var}\,B^i$ where ${\rm
Var}\,B^i$ is the total variation of the scalar process $B^i$. Obviously,
$$
|\dot B|_1=\sum_i|\dot B^i|=\sum_i
\left|\frac{dB^i}{d||B||}\right|=\sum_i \left|\frac{dB^i}{d{\rm
Var}\, B^i}\right|\frac{d{\rm Var}\,B^i}{d||B||}=\frac{d\sum_i{\rm
Var}\,B^i}{d||B||}.
$$
But  all norms in ${\bf R}^d$ are equivalent, i.e. $\tilde
\kappa^{-1}|.|\le |.|_1\le \tilde \kappa |.|$ for some strictly
positive constant $\tilde \kappa$. The same inequalities relate the corresponding total variation processes.
The claimed property follows from here with the constant $\kappa_1=\tilde \kappa^{-2}$.

Summarizing, we conclude that it is sufficient to check the
domination property for $E\tilde N_t$ with
\beq \label{expres}
\tilde N_t
:=t\wedge \tau + \int_0^{t\wedge \tau}I_{\{|c_s|\ge \kappa\}}
|c_s|ds +  ||B||_{t\wedge \tau}. \eeq
These  processes
$\tilde N=\tilde N^\pi$ have a transparent dependence on the control.  The idea of the concluding reasoning is very simple: on a certain
set of strictly positive probability, where one may neglect the
random fluctuations, either $\tau$ is ``large", or the total
variation of the control is ``large": one can accelerate
exit only by an intensive trading or consumption.

The formal arguments are as follows.  Using the
stochastic Cauchy formula
(\ref{dynam1}) and  the fact
 that $\cE_{0+}(Y^i)=\cE_0(Y^i)=1$, we get immediately
  that there
exist  a  number $t_0>0$   and a measurable set $\Gamma$ with
$P(\Gamma)>0$ on which
$$|V^{x,\pi}-x|\le
r/2 +2 (||B||+||C||) \quad \hbox{on}\ [0,t_0]
$$
whatever  is the control $\pi=(B,C)$. Of course, diminishing
$t_0$, we may assume without loss of generality that $\kappa
t_0\le r/8$. For any $t\le t_0$ we have on the set $\Gamma
\cap \{\tau \le t\}$ the inequality $||B||_\tau+||C||_{\tau}\ge
r/4$ and, hence,
$$
\tilde N_t\ge ||B||_\tau+||C||_{\tau}-\int_0^{\tau}I_{\{|c_s|<
\kappa\}} |c_s|ds\ge \frac {r}{4}-\kappa t_0\ge \kappa t_0\ge t_0\ge  t.
$$
On the set $\Gamma \cap \{\tau > t\}$ the inequality $\tilde N_t\ge
t$ is obvious. Thus, $E\tilde N_t\ge tP(\Gamma)$ on $[0,t_0]$ and the result is proven. \fdem

\section{Dynamic Programming Principle}



\noindent
The aim of this section is to establish the following two assertions which will serve to derive the HJB equation for the Bellman function. For the considered model, they constitute an analog of the classical Dynamic Programming Principle. The latter is usually written in the form of a single identity (see the remark at the end of the section), but for our purpose this  form, more precise, is needed.
\begin{lemm}
\label{upperb}
 Let ${\cal T}_f$  be the sets of
 finite  stopping times. Then
 \beq
 \label{wleq}
 W(x)\le \sup_{\pi\in \cA^x_a}\inf_{\tau\in {\cal T}_f}E\left (J^\pi_{\tau}+
e^{-\beta \tau}W(V^{x,\pi}_{\tau})I_{\{\tau< \theta\}} \right). \eeq
\end{lemm}
\begin{lemm}
\label{lowb} Suppose  that $W$ is continuous on ${\rm int}\,
K$. Then
 for any $\tau\in {\cal T}_f$
 \beq
 \label{wgeq}
 W(x)\ge \sup_{\pi\in \cA^x_a}E\left (J^\pi_{\tau}+
e^{-\beta \tau}W(V^{x,\pi}_{\tau})I_{\{\tau< \theta\}} \right). \eeq
\end{lemm}

We work on the canonical filtered space of c\`adl\`ag  functions
equipped with the measure $P$ which is the distribution of the driving L\'evy process. The generic point
$\omega=\omega_.$ of this space is a  $d$-dimensional c\`adl\`ag function on ${\bf
R}_+$, zero at the origin. Let $\cF^\circ_t := \sigma \{\omega_s,
\ s\le t\}$ and $\cF_t:= \cap_{\e>0}\cF^\circ_{t+\e}$. We add the
superscript $P$ to denote $\sigma$-algebras augmented by all
$P$-null sets from $\Omega$. Recall that $\cF^{\circ,P}_t$
coincides with $\cF^P_t$ (this assertion follows easily from the
predictable representation theorem).
The Skorohod metric makes $\Omega$ a Polish space and its Borel $\sigma$-algebra coincides with $\cF_\infty$, for details see \cite{JS:87}.

Since elements of $\Omega$ are paths,  we can  define such
operators as the stopping $\omega_.\mapsto \omega^s_.$, $s\ge 0$,
where $\omega^s_.=\omega_{s\wedge .}$ and the translation
$\omega_.\mapsto \omega_{s+.}-\omega_s$. Taking Doob's theorem
into account, one can describe $\cF^\circ_s$-measurable random
variables as those of the form $g(w_.)=g(w^s_.)$ where $g$ is a
measurable function on $\Omega$.

We define also the ``concatenation" operator as the measurable
mapping
$$
g:{\bf R}_+\times \Omega\times \Omega\to \Omega
$$
with $g_t(s,\omega_.,\tilde\omega_.)=\omega_t I_{[0,s[}(t)+(\tilde
\omega_{t-s}+\omega_s)I_{[s,\infty[}(t)$.

Notice that
$$
g_t(s,\omega^s_.,\omega_{.+s}-\omega_s)=\omega_t .
$$
Thus, $\pi (\omega)=\pi (g(s,\omega^s_.,\omega_{.+s}-\omega_s))$.

Let $\pi$ be a fixed strategy from $\cA^x_a$ and let
$\theta=\theta^{x,\pi}$ be the exit time from ${\rm int}\, K$ for the
process $V^{x,\pi}$.

\smallskip

Recall  the following general fact on regular conditional distributions.

Let $\xi $ and $\eta$ be two random variables taking values in
Polish spaces $X$ and  $Y$ equipped with their Borel
$\sigma$-algebras ${\cal X}$ and $\cY$.
 Then
$\xi$ admits a regular conditional distribution given $\eta=y$
which we shall denote by $p_{\xi |\eta }(\Gamma,y)$. This means that $p_{\xi |\eta }(.,y)$ is a probability measure on $\cX$, $p_{\xi |\eta }(\Gamma,.)$ is a $\cY$-measurable function, and
$$
E(f(\xi,\eta)|\eta)=\int f(x,y)p_{\xi |\eta}(dx,y)\bigg |_{y=\eta}
\quad (a.s.)
$$
for any $\cX\times\cY$-measurable function $f(x,y)\ge 0$.

We shall apply the above relation to the random variables
$\xi=(\omega_{.+\tau}-\omega_{\tau})$
 and
$\eta=(\tau,\omega^\tau)$. It is well-known that the L\'evy process starts afresh at stopping times, i.e. the  measure $P(.)$ itself (not depending on $y$) is the regular
conditional distribution $p_{\xi |\eta }(.,y)$.

At last, for fixed $s$ and $w^s$,   the shifted control $\pi_{.+s}
(g(s,\omega^s_.,\tilde \omega.))$  is admissible for the
initial condition $V^{x,\pi}_{s}(\omega)$ when $s\le \theta(\omega)$. Here we denote by
$\tilde \omega.$ a generic point of the canonical space.

\noindent
{\sl Proof of Lemma \ref{upperb}.} For arbitrary $\pi \in \cA^x_a$ and ${\cal T}_f$ we
have that
\bean
EJ^{\pi}_\infty&=& EJ^{\pi}_\tau + E e^{-\beta\tau}I_{\{\tau< \theta\}}\int_0^\infty
e^{-\beta r} U(c_{r+\tau})dr \\
&=&EJ^{\pi}_\tau + E e^{-\beta\tau}I_{\{\tau< \theta\}}E\Big(\int_0^\infty e^{-\beta r} U(c_{r+\tau})dr\Big |(\tau,\omega^\tau)\Big)
.
\eean
According to the above  discussion we can rewrite the second term
of the right-hand side as
$$
Ee^{-\beta \tau}I_{\{\tau< \theta\}}\int \left ( \int_0^\infty e^{-\beta r
}U(c_{r+\tau}(g(\tau,\omega^{\tau},\tilde \omega)))dr
\right)P(d\tilde \omega)
$$
and dominate it  by $ Ee^{-\beta \tau}I_{\{\tau< \theta\}}W(V^{x,\pi}_{\tau})$. Thus,
$$
EJ^{\pi}_\infty\le EJ^{\pi}_\tau + E
e^{-\beta\tau}I_{\{\tau< \theta\}}W(V^{x,\pi}_{\tau}).
$$
This bound leads directly to the  announced inequality. \fdem

\smallskip


\noindent
{\sl Proof of Lemma \ref{lowb}.} Fix $\e>0$. By hypothesis, the function $W$ is
continuous on ${\rm int}\,K$. For each
 $x \in {\rm int}\, K$ we
can find an open ball $\cO_r(x)=x+\cO_r(0)$ with $r=r(\e,x)<\e$
contained in   the open set $\{y\in {\rm int}\,K:\
|W(y)-W(x)|<\e\}$. Moreover, we can find
 a
smaller ball $\cO_{\tilde r}(x)$ contained in the set $y(x)+K$
with some
 $y(x)\in \cO_r(x)$. Indeed, take an arbitrary $x_0\in {\rm int}\,K$. Then, for some $\delta>0$, the ball $x_0+\cO_\delta (0)\subset
 K$.
 Since $K$ is a cone,  $\lambda x_0+\cO_{\lambda\delta} (0)\subset
 K$ for every $\lambda>0$ and this inclusion implies that
$$ x+\cO_{\lambda
\delta}(0)\subset x-\lambda x_0+K
$$
Clearly, the requirement is met for
$y(x)=x-\lambda x_0$ and $\tilde r=\lambda \delta $ when
$\lambda|x_0|<r$ and $\lambda \delta <r$. The family of sets
$\cO_{\tilde r(x)/2}(x)$, $x\in {\rm int}\, K$, is an open covering
of ${\rm int}\, K$. But any open covering of a separable metric
space  contains a countable subcovering (this is the Lindel\"of
property; in our case, where  ${\rm int}\, K$ is a countable union
of compacts, it is obvious). Take a countable subcovering indexed
by points $x_n$. For notational simplicity, we shall denote the open balls 
$\cO_{\tilde r(x_n)/2}(x_n)$
by
$\cO_n$ and $y(x_n)$ by $y_n$. 

Let $\pi^n=(B^n,C^n)\in \cA^{y_n}_a$ be an $\e$-optimal strategy
for  the initial point $y_n$, i.e. such that
$$
EJ^{\pi_n}_\infty\ge W(y_n)-\e.
$$
Let $\pi\in \cA^x_a$
 be an arbitrary strategy. Put 
$$
 \rho:=\inf \{j\ge 1\colon\ V^{x,\pi}_{\tau}\in \cO_j\},  
 $$
 Let us introduce  the strategy 
 $$
 \pi':=\pi I_{[0,\tau]}+(0,0)I_{]\tau,\infty[} 
 $$ 
 and the predictable stopping times $\tau_k:=\tau+1/k$. 
 Finally, put  
\bean
\tilde \pi:&=&\pi I_{[0,\tau]}
+\sum_{n=1}^\infty
[(y_n-V^{x,\pi'}_{\tau_k},0)+\bar \pi^{n,k}
]I_{[\tau_k,\infty[}I_{\{\rho=n\}}
I_{\{V^{x,\pi'}_{\tau_k}-y_n\in K\}}
I_{\{\tau_k< \theta\}}
\eean
where $\bar \pi^{n,k}$ is the translation of the strategy  $\pi^n$:
namely, for
 a point  $\omega_.$ with  $\tau (\omega_.)=s<\infty$ we have
  $$
  \bar \pi^{n,k}_t(\omega_.):=\pi^n_{t-s-1/k} (\omega_{.+s+1/k}-\omega_{s+1/k}), \quad t\ge s_1/k. 
  $$
In other words, the strategy $\tilde \pi$ coincides with $\pi$ on
$[0,\tau[$,  is zero on the interval $[\tau,\tau_k[$ 
 and with the shift of $\pi^n$ on  $[\tau_k,\infty[$ when $V^{x,\pi}_{\tau}$
 is in  $\cO_n$ and $V^{x,\pi}_{\tau_k}-y_n\in K$; the correction term guarantees that
in the latter case the trajectory
 of the control system corresponding to the control $\tilde \pi$
 passes at time $\tau_k$ through the point $y_n$.
 
 One can check that $\tilde \pi I_{[0,\theta^{x,\tilde \pi}]}\in \cA_a^x$. 

 Now, using the same considerations as in the previous lemma, we
 have:
 \bean
 W(x)\ge EJ^{\tilde \pi}_{\infty}&=&E J^{\pi}_\tau+\sum_{n=1}^\infty
 E I_{\{\rho=n\}}I_{\{\tau< \theta\}}I_{\{V^{x,\pi'}_{\tau_k}-y_n\in K\}}\int _{\tau_k}^\infty
 e^{-\beta s}U(\bar c^n_s)ds\\
 & \ge &E J^{\pi}_\tau+\sum_{n=1}^\infty E I_{\{\rho=n\}}
 I_{\{\tau< \theta\}}e^{-\beta \tau }
 (W(y_n)-\e)\\
 &\ge& E J^{ \pi}_{\tau}+ Ee^{-\beta \tau }W(V^{x,\pi}_{\tau})I_{\{\tau< \theta\}} - 2\e.
 \eean
Since $\pi$ and $\e$ are arbitrary, the result follows.  \fdem

\smallskip
\noindent{\bf Remark.} The previous lemmas imply that for any $\tau\in {\cal T}_f$ the following identity holds:
$$
 W(x)= \sup_{\pi\in \cA^x_a}E\left (J^\pi_{\tau}+
e^{-\beta \tau}W(V^{x,\pi}_{\tau})I_{\{\tau< \theta\}} \right).
$$
It can be considered as a form of the dynamic programming
principle but, seemingly, it is not sufficient for our derivation of the
HJB equation.


\section{The Bellman Function and the HJB Equation}
\begin{theo}
\label{conditional}
 Assume
that the Bellman function  $W$ is in $C(K)$. Then $W$ is a
viscosity solution of\, (\ref{vi}).
\end{theo}
{\sl Proof.} The claim follows  from the two lemmas below. \fdem

\begin{lemm} If $W$ is in $C({\rm int}\, K)$ 
then  $W\ge 0$ is a viscosity supersolution of (\ref{vi}).
\end{lemm}
{\sl Proof.} Let $x\in$ int $K$ and let $\phi\in C^1(K)\cap C^2(x)$ be a function such that $\phi(x)=W(x)$ and $W\ge \phi$ on $K$. 
  
Fix an arbitrary point $m\in K$. Let $\e>0$ be sufficiently small to
guarantee that $x-\e m \in \cO_r(x)$. The function $W$ is increasing with
respect to the partial ordering generated by $K$. Thus,
$$
\phi(x)=W(x)\ge W(x-\e m)\ge \phi (x-\e m).
$$
Taking a limit as $\e\to 0$ in the inequality 
$\e^{-1}(\phi (x-\e m)-\phi (x))\le 0$,
we obtain that
$-m\phi'(x)\le 0$ and, hence, $\Sigma_G(\phi'(x))\le 0$.

Take now $\pi$ with $B_t=0$ and  $c_t=c\in \cC$ for all $t$. Let $\tau_r=\tau^\pi_r \le \theta$ be
the exit time of the  process $V=V^{x,\pi}$ from the
ball $\cO_r(x)$; obviously, $\tau_r\le \theta$.  The properties of the test function and the inequality
(\ref{wgeq}) imply that
\bean
\phi(x)=W(x)&\ge& E\left (J^\pi_{t\wedge\tau_r}+ e^{-\beta
(t\wedge\tau_r)}W(V_{t\wedge\tau_r})I_{\{t\wedge\tau_r< \theta\}} \right)\\&\ge& E\left (J^\pi_{t\wedge\tau_r}+ e^{-\beta
(t\wedge\tau_r)}\phi(V_{t\wedge\tau_r})I_{\{t\wedge\tau_r< \theta\}} \right).
\eean
We get from here, using the
It\^o formula (\ref{Ito}), that
\bean
0&\ge& E\left(\int_0^{t\wedge \tau_r}e^{-\beta s}U(c_s)ds+
e^{-\beta(t\wedge \tau_r)}\phi(V_{t\wedge\tau_r})I_{\{t\wedge\tau_r< \theta\}} \right )-\phi (x)\\
&\ge& EI_{\{t\wedge\tau_r< \theta\}}\int_0^{t\wedge \tau_r}e^{-\beta
s}[\cL_0\phi(V_s)-c\phi'(V_s)+U(c)
]ds\\
&\ge& \min_{y\in \bar \cO_r(x)} [\cL_0\phi(y)-c \phi'(y)+U(c)]
EI_{\{t\wedge\tau_r< \theta\}}\left [\frac{1}{\beta}\left(1-e^{-\beta(t\wedge
\tau_r)}\right)\right].
\eean
Dividing the resulting inequality by
$t$ and taking successively the limits as  $t$ and $r$ converge to
zero we infer that $\cL_0\phi(x)-c \phi'(x)+U(c)\le 0$. Maximizing
over $c\in \cC$ yields the bound  $\cL_0\phi(x)+U^*(\phi'(x))\le
0$ and, therefore, $W$ is a supersolution of the HJB equation.
 \fdem

\begin{lemm} If (\ref{wleq}) holds
then  $W\ge 0$ is a viscosity subsolution of (\ref{vi}).
\end{lemm}
{\sl Proof.} Let $x\in {\rm int}\, K$   and let $\phi\in C^1(K)\cap C^2(x)$ be a function such that $\phi(x)=W(x)$ and $W\le \phi$
on $K$.  Suppose  that the subsolution inequality for $\phi$ fails
at $x$. Thus, there exists $\e>0$ such that $\cL\phi\le -\e$ on
some ball $\bar \cO_r(x)\subset {\rm int}\, K$. By virtue of Lemma
\ref{svi} (applied to the function $\phi$) there are $t_0>0$ and
$\eta>0$ such that on the interval $]0,t_0]$ for any strategy
$\pi\in \cA_a^x$
$$
E\left (J^\pi_{t\wedge \tau^\pi}+ e^{-\beta
\tau^\pi}\phi(V^{x,\pi}_{t\wedge \tau^\pi})I_{\{t\wedge \tau^\pi< \theta\}} \right)\le
\phi(x)-\eta t,
$$
where $\tau^\pi$ is the exit time of the process $V^{x,\pi}$ from
the ball $ \cO_r(x)$. Fix arbitrary  $t\in ]0,t_0]$. By the second claim
of Lemma \ref{upperb} there exists $\pi\in \cA_a^x$ such that
$$
\label{supremum} W(x)\le E\left (J^\pi_{t\wedge \tau}+ e^{-\beta
\tau}W(V^{x,\pi}_{t\wedge \tau})I_{\{t\wedge \tau< \theta\}} \right)+\frac 12 \eta t,
 $$
 for every  stopping time $\tau$, in particular for $\tau^\pi$.

Using the inequality $W\le \phi$ and applying Lemma \ref{svi} we
obtain from the above relations that $W(x)\le \phi (x)-(1/2)\eta
t$. This is a contradiction because at the point $x$ the values of
$W$ and $\phi$ are the same. \fdem

\section{Uniqueness Theorem}
Before formulating the uniqueness theorem we recall the Ishii lemma.

\begin{lemm}
\label{Ishii} Let $v$ and $\tilde v$ be two continuous functions
on an open subset $\cO\subseteq {\bf R}^d$. Consider the function $\Delta (x,y):= v(x)-\tilde
v(y)-\frac 12 n|x-y|^2$ with $n>0$. Suppose that $\Delta$ attains
a local maximum at $(\widehat x,\widehat y)$. Then there are
symmetric matrices $X$ and $Y$ such that
$$
(n(\widehat x-\widehat y),X)\in \bar J^+v(\widehat x), \qquad
(n(\widehat x-\widehat y),Y)\in \bar J^-\tilde v(\widehat y),
$$
and \beq \label{doublematrix} \left (\begin{array}{cc} X&\ 0\\
0&-Y
\end{array}
\right) \le 3n \left (\begin{array}{cc} \ I&-I\\ -I&\ I
\end{array}
\right). \eeq
\end{lemm}

In this statement $I$ is the identity matrix and  $\bar J^+v$
and $\bar J^-\tilde v$ are values of the set-valued mappings whose
graphs  are closures of graphs of the set-value mappings $J^+v$ and $J^-\tilde v$, respectively.

The inequality (\ref{doublematrix}) means that for any vectors $x$ and $y$ from ${\bf R}^d$
\beq
(Xx,y)-(Yy,y)\le 3n|x-y|^2.
\eeq

Of course, if $v$ is smooth, the claim follows directly from the
necessary conditions of a local maximum (with $X=v''(\widehat x)$,
$Y=\tilde v''(\widehat y)$ and the constant $1$ instead of $3$ in
inequality (\ref{doublematrix})).

The inequality (\ref{doublematrix}) implies the  bound
\beq\label{doublematrix1}
{\rm tr}\,(A(x)X -A(y)Y)\le 3n|A^{1/2}|^2|x-y|^2
\eeq
which will be used in the sequel (for the proof see, e.g., Section 4.2 in \cite{K-Saf}).

The following concept plays a crucial role in the proof of the
purely analytic result on the uniqueness of the viscosity solution
which we establish by a classical method of doubling variables
using the Ishii lemma.
\smallskip

\noindent
 {\bf Definition.} We say that a positive function $\ell\in C_1(K)\cap C^2({\rm
int}\, K $) is the {\it Lyapunov
 function}\index{Lyapunov
 function} if the following properties are satisfied:

1)
  $\ell'(x)\in {\rm int}\,K^*$ and $\cL_0 \ell(x)\le 0$ for all  $x\in {\rm int}\,
K $,

2) $\ell(x)\to \infty$ as $|x|\to \infty$.

\smallskip In other words, $\ell$ is a classical strict
supersolution of the {\bf truncated} equation (without
the term $U^*$), continuous up to the boundary, and increasing to infinity at infinity.

\begin{theo}
\label{unique}
 Suppose that there exists a  Lyapunov function
$\ell$ and the L\'evy measure $\Pi$  is such that 
$$
\Pi(z:\ \hat x+D_{\hat x}z\in \partial K\})=0 \qquad  \forall \hat x\in {\rm int\,}K.
$$
Then
the Dirichlet problem (\ref{vi}), (\ref{icond})  has at most one
viscosity solution in the class of continuous functions satisfying
the growth condition \beq W(x)/\ell (x)\to 0, \quad |x|\to \infty.
\eeq
\end{theo}

\noindent {\sl Proof.} Let $W$ and $\tilde W$ be two  viscosity
solutions
 of (\ref{vi}) coinciding on the boundary $\partial K$. Suppose that
$W(z)>\tilde W(z)$ for some $z\in K$. Take $\e>0$ such that
$$
W(z)-\tilde W(z)-2\e\ell(z)>0.
 $$
We introduce a family of continuous functions $\Delta_n:K\times K\to {\bf R}$
 by putting
$$
\Delta_n(x,y):=W(x)-\tilde W(y)-\frac 12 n|x-y|^2-\e
[\ell(x)+\ell(y)], \qquad n\ge 0.
$$
Note that $\Delta_n(x,x)=\Delta_0(x,x)$ for all $x\in K$ and
$\Delta_0(x,x)\le 0$ when $x\in \partial K$. From the assumption
that the function $\ell$ has a higher growth rate than $W$ we deduce
that $\Delta_n(x,y)\to -\infty$ as $|x|+|y|\to \infty$. It follows
that the level sets $\{\Delta_n\ge a\}$ are compacts and the
function $\Delta_n$ attains its maximum on a compact subset of $K\times K$ which does not depend on $n$. That is, there exists
$(x_n,y_n)\in K\times K$ such that
$$
\Delta_n(x_n,y_n)=\bar \Delta_n:=\sup_{(x,y)\in K\times
K}\Delta_n(x,y)\ge \bar \Delta :=\sup_{x\in K}\Delta_0(x,x)>0.
$$
All $(x_n,y_n)$ belong to the compact set $\{(x,y):\
\Delta_0(x,y)\ge 0\}$. It follows that the sequence $n|x_n-y_n|^2$
is  bounded. We continue to argue (without introducing new
notations) with a subsequence along which $(x_n,y_n)$ converge to
some limit $(\widehat x,\widehat x)$. Necessarily,
$n|x_n-y_n|^2\to 0$ (otherwise we would have
 $\Delta_0 (\widehat x,\widehat x)>\bar \Delta$). It is easily seen that
$\bar \Delta_n \to \Delta_0(\widehat x,\widehat x)=\bar \Delta$.
Thus,
 $\widehat x$ is an interior  point  of
$K$ and so are  $x_n$ and $y_n$ for sufficiently large $n$. Let $j\ge 1$ be the number of 
nonzero components  $\widehat x$. Without loss of generality we assume that 
$\widehat x^{j+1},\dots,\widehat xd=0$ and for sufficiently large $n$ the first $j$ components of $x_n$ and $y_n$ are strictly positive.    

\smallskip

By virtue of the Ishii lemma applied to the functions $v:=W-\e\ell
$ and $\tilde v:=\tilde W+\e\ell $ at the point $(x_n,y_n)$ there
exist matrices  $X^n$ and $Y^n$ satisfying (\ref{doublematrix})  such
that
\beq
\label{natural}
(n(x_n-y_n),X^n)\in \bar J^+v(x_n), \qquad (n(x_n-y_n),Y^n)\in \bar
J^-\tilde v(y_n).
\eeq
Suppose for a moment that
\beq
\label{unnatural}
(n(x_n-y_n),X^n)\in  J^+v(x_n), \qquad (n(x_n-y_n),Y^n)\in
J^-\tilde v(y_n).
\eeq
Using the notations  $ p_n:=n(x_n-y_n)+\e \ell'(x_n)$,
$q_n:=n(x_n-y_n)-\e \ell'(y_n)$ and putting  $X_n:=X^n+\e \ell''(x_n)$,
$Y_n:=Y^n-\e \ell''(y_n)$, we may rewrite the last relations in the
following equivalent form: \beq (p_n,X_n)\in  J^+W(x_n),
\qquad (q_n,Y_n)\in  J^-\tilde W(y_n). \eeq 

Since $W$  is a viscosity subsolution, by virtue of Lemma \ref{jet} there exists a function $f_n\in C_1(K)\cap C^2(x_n)$  such that 
$f'_n(x_n)=p_n$, $f''_n(x_n)=X_n$, $f_n(x_n)=W(x_n)$, and 
$W\le f_n\le W+1/n$ on $K$. 
Since $\tilde W$  is a viscosity supersolution we conclude in the same way that  there exists a function $\tilde f_n\in C_1(K)\cap C^2(y_n)$ such that 
$\tilde f'_n(y_n)=q_n$, $\tilde f''_n(y_n)=Y_n$, $\tilde f_n(y_n)=\tilde W(y_n)$,
and $\tilde W-1/n \le \tilde f_n\le \tilde  W$ on $K$. 
      To deal with the nonlocal integral operator we take $f_n$ and $\tilde f_n$ having the structure given in (\ref{cstrf})  with an appropriate choice of the cylindrical cutoff functions. We discuss details of this choice  later.  

By definitions of sub- and supersolutions we have that 
$$
F(X_n,p_n,\cI(f_n,x_n), W(x_n),x_n)\ge 0\ge F(Y_n,q_n,\cI(\tilde f_n,y_n),\tilde W (y_n),y_n).
$$
The second inequality implies that  $mq_n\le 0$ for each $m\in
G=(-K)\cap \partial  \cO_1(0)$. But for the Lyapunov function
$\ell'(x)\in {\rm int}\, K^*$ when $x\in {\rm int}\, K$ and,
therefore,
$$
m p_n= mq_n+\e m(\ell'(x_n)+\ell'(y_n))<0.
$$
Since $G$ is a compact, $\Sigma_G(p_n)<0$. It follows that
\bean
F_0(X_n,p_n,\cI(f_n,x_n), W(x_n),x_n)+U^*(p_n)&\ge& 0,\\
F_0(Y_n,q_n,\cI(\tilde f_n,y_n),\tilde W (y_n),y_n)+U^*(q_n)&\le& 0.
\eean
Recall that $U^*$ is  decreasing with respect to the partial
ordering generated by $\cC^*$ hence also by $K^*$. Thus,
$U^*(p_n)\le U^*(q_n)$ and we obtain the inequality
$$
b_n:=F_0(X_n,p_n,\cI(f_n,x_n), W(x_n),x_n)- F_0(Y_n,q_n,\cI(\tilde f_n,y_n),\tilde W (y_n),y_n)\ge
0.
$$
Clearly,
 \bean b_n&=&\frac 12 \sum_{i,j=1}^d(a^{ij}x^i_{n}x^j_n
X^n_{ij} - a^{ij}y^i_{n}y^j_n Y^n_{ij})+n \sum_{i=1}^d \mu^i
(x^i_n-y^i_n)^2\\&& - \frac 12\beta n|x_n-y_n|^2
-\beta\Delta_n(x_n,y_n) +\cI(f_n-\e \ell,x_n)-\cI(\tilde f_n+\e\ell,y_n)\\
&&+\e(\cL_0\ell (x_n)+\cL_0\ell (y_n)).
\eean
By virtue of (\ref{doublematrix1}) the first term in the right-hand side is
dominated by a constant multiplied by $n|x_n-y_n|^2$; a similar bound for
the second sum is obvious; the last term is negative according to
the definition of the Lyapunov function. To complete the proof, it is sufficient to show that
\beq\label{key-int}
\limsup_n\,(\cI(f_n-\e \ell,x_n)-\cI(\tilde f_n+\e\ell,y_n))\le 0.
\eeq
Indeed, with this we have that  $\limsup
b_n\le -\beta \bar \Delta<0 $, i.e.  a contradiction arising
from the assumption $W(z)>\tilde W(z)$.

In general,  we cannot guarantee that (\ref{key-int}) holds for arbitrary test functions 
$f_n$ and $\tilde f_n$.
That is why we choose them in accordance with the expressions given by Lemma 
\ref{jet} with $\alpha=1/n$, i.e. 
\bean
\label{appro1}
f_n(x_n+h)&:=&f^0_n(x_n+h)\xi^L_{a_n,a'_n}(h)+W(x_n+h)\big(1-\xi^L_{a_n,a'_n}(h)\big), \quad x_n+h\in K,\qquad \\
\label{appro2}
\tilde f_n(y_n+h)&:=&\tilde f^0_n(y_n+h)\tilde \xi^L_{\tilde a_n,\tilde a'_n}(h)+\tilde W(y_n+h)\big(1-\tilde \xi^L_{\tilde a_n,\tilde a'_n}(h)\big),\quad y_n+h\in K,\qquad 
\eean
where the linear space $L:={\rm Im}\, D_{\widehat x}=\{x\in {\bf R}^d\colon\ x^{j+1}=0, \dots,x^d=0 \}$,  
\bean
\label{appro3}
f^0_n(x_n+h)&:=&\big[\big(W(x_n)+Q_{p_n,X_n}(h)+r_n(h)\big)\vee W_n(x_n+h)\big]\wedge \big[W(x_n+h)+1/n\big],\\
\label{appro4}
\tilde f^0_n(y_n+h)&:=&\big[\big(\tilde W(y_n)+Q_{q_n,Y_n}(h)-\tilde r_n(h)\big)\wedge \tilde W_n(y_n+h)\big]\vee \big[\tilde W(y_n+h)-1/n\big].  
\eean

 Let $\delta:=(1/2)R/(|\widehat x|+R)$ where 
$R=d(\widehat x, \partial K)$ is the distance of the point $\widehat x$ from the boundary $\partial K$.
Then $I(z,x_n)=1$ and $I(z,y_n)=1$ for $z\in \cO_\delta(0)$ when $n$ is sufficiently large. Indeed, when 
$|x_n-\widehat x|\le R/2$ and $|y_n-\widehat x|\le R/2$ we have  
that 
$$
|x_n +D_{x_n}z -\widehat x|\le R/2+|x_n||z|\le R/2+(|\widehat x|+R/2)|z|< R 
$$
when $|z|\le \delta$, and the similar estimate holds for $y_n$. 

We have:   
$$\cI(f_n-\e \ell,x_n)-\cI(\tilde f_n+\e\ell,y_n)=\int_{\{|z|\le \delta\}}
H_n(z)\Pi(dz)+\int_{\{|z|> \delta\}}
H_n(z)\Pi(dz)
$$
where $H_n(z):=F_n(z)-\tilde F_n(z)$ with 
\bean
F_n(z)&:=&(f_n-\e\ell)(x_n+D_{
x_n} z)I(z,x_n)-(W-\e\ell)(x_n)-(f_n'-\e\ell')(x_n)D_{x_n} z\\
&=&(f_n-\e\ell)(x_n+D_{
x_n} z)I(z,x_n)-(W-\e\ell)(x_n)-n(x_n-y_n)D_{
x_n} z
,\\
\tilde F_n(z)&:=&(\tilde f_n+\e\ell)(y_n+D_{y_n} z)I(z,y_n)-(\tilde W+\e\ell)(y_n)-(\tilde f'_n+\e\ell')(y_n)D_{y_n} z\\
&=&(\tilde f_n+\e\ell)(y_n+D_{
y_n} z)I(z,y_n)-(\tilde W+\e\ell)(y_n)-n(x_n-y_n)D_{
y_n} z.
\eean

Define also the functions
\bean
G_n(z)&:=&(W-\e\ell)(x_n+D_{
x_n} z)I(z,x_n)-(W-\e\ell)(x_n)-n(x_n-y_n)D_{
x_n} z
,\\
\tilde G_n(z)&:=&(\tilde W+\e\ell)(y_n+D_{
y_n} z)I(z,y_n)-(\tilde W+\e\ell)(y_n)-n(x_n-y_n)D_{
y_n} z.
\eean

 
Put $x'_n:=P_Lx_n$,  $y'_n:=P_Ly_n$.  
Considering only sufficiently large 
$n$,  we may assume without loss of generality that $(1/2)\widehat x^i 
\le x_n^i,\,  y_n^i \le \widehat x^i+1$ for non-zero coordinates of $\widehat x$ and, therefore, the norms of the restrictions of the diagonal operators   $D_{x'_n}$, $D_{y'_n}$ on $L$ and their inverses are bounded (even uniformly in $n$). Therefore, we may apply Lemma \ref{Mollif=} if $\widehat x$ has no zero components or its extension  given by the accompanying  remark  and  argue further supposing that 
$$
\xi^L_{a_n,a'_n}(D_{x_n}z)=\xi^L_{a_n,a'_n}(D_{x'_n}z)=\xi^L_{\tilde a_n,\tilde a'_n}(D_{y'_n}z)=\xi^L_{\tilde a_n,\tilde a'_n}(D_{y_n}z) \quad \forall\, z\in {\bf R}^d. 
$$
According to our choice of $\delta$ for $z\in \cO_\delta(0)$ we have  $I(z,x_n)=I(z,y_n)=1$ for sufficiently large $n$ and, as a consequence, the following easily verified identity:   
$$
G_n(z)-\tilde G_n(z)=\Delta_n(x_n+D_{x_n}z,y_n+D_{y_n}z) - \Delta_n(x_n,y_n)
+\frac 12 n|(D_{x_n}-D_{y_n})z|^2. 
$$
Recalling that  the function $\Delta_n(x,y)$ attains its maximum at $(x_n,y_n)$ we get from here the bound
$$
G_n(z)-\tilde G_n(z)\le \frac 12 n|x_n-y_n|^2|z|^2. 
$$
Thus, on the set $\{z:\ |D_{x'_n}z|\ge a'_n \}\cap \cO_\delta(0)$ we have that 
\beq
H_n(z)=G_n(z)-\tilde G_n(z)\le \frac 12 n|x_n-y_n|^2|z|^2.
\eeq
On the set $\{z:\ |D_{x'_n}z|\le a_n \}\cap \cO_\delta(0)$ 
\bean
F_n(z)&=&\frac 12(X^nD_{x_n}z,D_{x_n}z)+r_n(D_{x_n}z)\\&&
+\e\Big[\ell(x_n)+\ell '(x_n)D_{x_n}z+\frac 12(\ell''(x_n)D_{x_n}z,D_{x_n}z)-\ell(x_n+D_{x_n}z)\Big],\\ 
\tilde F_n(z)&=&\frac 12(Y^nD_{y_n}z,D_{y_n}z)+\tilde r_n(D_{y_n}z)\\&&
+\e\Big[\ell(y_n)+\ell '(y_n)D_{y_n}z+\frac 12(\ell''(y_n)D_{y_n}z,D_{y_n}z)-\ell(y_n+D_{y_n}z)\Big]. 
\eean
By the Ishii lemma
$$
(X^nD_{x_n}z,D_{x_n}z)-(Y^nD_{y_n}z,D_{y_n}z)\le 3n|D_{x_n}z-D_{y_n}z|^2\le 3n|x_n-y_n|^2|z|^2. 
$$ 
We take $a_n$ small enough to ensure that 
$$
\int_{\{|D_{x'_n}z|\le a_n\}} (r(D_{x_n}z)+r(D_{x_n}z)\Pi(dz)\le 1/n. 
$$
On the set $\{z:\ |D_{x'_n}z|\ge a_n \}\cap \cO_\delta(0)$ the sequence 
of functions $|H_n|$ is bounded by a constant. 
We choose
$a'_n$ sufficiently close to $a_n$ to guarantee that 
$$
\int _{\{a_n\le |D_{x'_n}z|\le a'_n\}}|H_n(z)|\Pi(dz)\le \frac 1n. 
$$
The expressions in parentheses in the  above formulae for $F_n(z)$ and 
$\tilde F_n(z)$ (residual terms in the Taylor formula for the smooth function $\ell$) are bounded by a constant  times $|z|^2$ and converge to the same limit as $n\to \infty$.  

Summarizing the above facts we conclude that 
$$
\limsup_n \int_{\{|z|\le \delta\}}
H_n(z)\Pi(dz)=0.  
$$

Since the continuous function $W$ and $\ell$ are of sublinear growth and the sequences $x_n$ and $n(x_n-y_n)$ are converging (hence bounded), the absolute value of $F_n$ is dominated by a function $c (1+|z|)$. The arguments for $-\tilde F_n(z)$ are similar. So, the function $H_n$ is dominated by a function of sublinear growth.

 Put $\mathcal{Z}:=\{z\,:\, \hat x+D_{\hat x}z\in \partial K\}$. By assumption, $\mathcal{Z}$ is $\Pi$-null. If $z\notin \mathcal{Z}$ and $n$ is large enough, we  have  $I(x_n,z)=I(y_n,z)=I(\hat x,z)$.  It follows that 
 \bean H_n(z)&=&(G_n(z)-\tilde G_n(z))I(\hat x,z)-n(x_n-y_n)(D_{x_n}-D_{y_n})z(1-I(\hat x,z))\\
&&-(W-f_n)(x_n+D_{x_n}z)I(\hat x,z)+(\tilde W-\tilde f_n)(y_n+D_{y_n}z)I(\hat x,z)
\\
&&+\big (\tilde W+\e\ell)(y_n)-(W-\e\ell)(x_n)\big)(1-I(\hat x,z))\\
&\le &
 \left(\Delta_n (x_n+D_{x_n}z, y_n+D_{y_n}z)-\Delta_n(x_n, y_n)\right)I(\hat x,z)\\
&+&\left(\frac 12 n|x_n-y_n+(D_{x_n}-D_{y_n})z|^2-\frac 12 n|x_n-y_n|^2\right)I(\hat x,z)\\&-&n(x_n-y_n)(D_{x_n}-D_{y_n})z+2/n\\
&&+\big (\tilde W+\e\ell)(y_n)-(W-\e\ell)(x_n)\big)(1-I(\hat x,z)).
\eean 
 
The first term in the right-hand side is negative as $\Delta_n(x,y)$ attains its maximum  at $(x_n,y_n)$. We only need to consider the second term when $I(\hat x,z)\ne 0$. In this case, we may combine it with the third term and conclude that the sum is less than $(1/2)n|x_n-y_n|^2\to 0$. The last term converges to zero because of continuity. Thus, 
$$
\limsup_n \int_{\{|z|> \delta\}}
H_n(z)\Pi(dz)=0.  
$$

However, our reasoning  is based on the assumption (\ref{unnatural}) while we know only
(\ref{natural}). Fortunately, using the definitions of $\bar J^+v(x_n)$ and $\bar J^-\tilde v(y_n)$ we can replace
the objects $x_n$, $y_n$, $X_n$, $Y_n$ by  their approximations 
$\widehat x_n$, $\widehat y_n$, $\widehat X_n$, $\widehat Y_n$ approaching rapidly the initial ones and for those  (\ref{unnatural}) hold. Repeating the arguments and controlling the approximation errors, we get the same contradiction.
 \fdem

\smallskip
\noindent {\bf Remark 1.}
Note that the definition of the Lyapunov function
does not depend on $U$  and hence the uniqueness holds for any utility function $U$
for which $U^*$ is decreasing with respect to the partial ordering
induced by $K^*$. However,  to apply the uniqueness theorem one needs to determine the growth rate of  $W$ and provide a Lyapunov function with a faster growth.

\section{Existence of Lyapunov Functions and Classical
Supersolutions}
In this section we extend results of \cite{K-Kl} on the existence of  Lyapunov functions and classical supersolutions   to the considered case of nonlocal operators. 

\smallskip
\noindent
{\bf Construction of Lyapunov functions.}
\smallskip

Let $u\in C({\bf R}_+)\cap C^2({\bf R}_+\setminus \{0\})$ be an
increasing strictly concave function with $u(0)=0$ and
$u(\infty)=\infty$. Introduce the function $R:=-u'^2/(u''u)$.
Assume that $\bar R :=\sup_{z>0} R(z)<\infty$.

For $p\in K^*$ we define on $K$ the positive function $f\in C_1(K)\cap C^2({\rm
int}\, K$ by putting  
$f(x)=f_p(x):=u(px)$. If $y\in K$, then
$yf'(x)=(py)u'(px)\ge 0$.

If $p\in {\rm int}\,K^*$, then for any $x,y\in K\setminus \{0\}$  we have the strict inequality $yf'(x)> 0$ implying that  $f'(x)\in {\rm int}\,K^*$. Thus, for
$p\in {\rm int}\,K^*$ the function $f$ is a Lyapunov function provided that the inequality $\cL_0f(x)\le 0$ is satisfied.   We show that under some mild conditions this inequality holds for sufficiently large $\beta$. 

Put  $\kappa_p:=0$, $\eta_p:=0$, if $\Pi=0$ and 
 $$
 \kappa_p:= \sup_{x\in {\rm int}\,K}\frac {u'(px)}{u(px)}|p||x|,  
\quad
\eta_p:=\kappa_p\int _{\{|z|>\kappa_p^{-1}\}} |z|\Pi(dz) 
$$
otherwise. 
Define also  
$$
\tilde \eta_p
:=\frac 12 \sup_{x\in {\rm int}\,K} \frac {\langle \mu
(x),p\rangle^2}{\langle A(x)p,p\rangle}I_{\{\langle A(x)p,p\rangle\neq 0\}}. 
 $$ 

Recall that  $A(x)$ is the matrix with $A^{ij}(x)=a^{ij}x^ix^j$
and the vector $\mu (x)$ has the components $\mu^ix^i$. 

Note that if   $ \kappa_p<\infty$, then  $\eta_p<\infty$ (as we assume  that $\int |z|^2\wedge | z|\Pi(dz)<\infty$). 

\smallskip
\noindent
{\sl Example.} Let $u(z):=z^\rho/\rho$ where $\rho\in ]0,1[$. Then $\bar R=R(z)=\rho/(1-\rho)$ and, for $p\in {\rm int}\,K^*$ 
$$
\kappa_p\le \rho \sup_{x\in K\setminus \{0\}} \frac{|p||x|}{px}<\infty 
$$
(strictly positive function $y\mapsto py$ on the compact $K\cap \{y\colon\ |y|=1\}$
attains its minimum). 

\begin{prop}
\label{Lyap} 
Let $p\in {\rm int}\,K$. If $\kappa_p<\infty$ and
$\beta\ge\tilde \eta_p \bar R +\eta_p+\max_i |\mu_i|\kappa_p$, then $f_p$ is a Lyapunov function.
\end{prop}
{\sl Proof.} Let $x\in  {\rm int}\,K$. Recall that 
$$
\cI(f,x): =\int \big[(f(x+D_xz)I_{{\rm int}\, K}(x+D_xz)-f(x))-D_xz f'(x)\big]\Pi(dz).
$$
If  $x+D_x z\in {\rm int}\,K$, then the integrand defining $\cI(f,x)$ has three nontrivial terms and we have by the Taylor formula (in which   $\vartheta \in [0,1]$) that
$$
(f(x+D_x  z)-f(x)-D_x  zf'(x))=\frac 12u''(px+\vartheta p\,D_x z)(p\,D_x  z )^2\le 0. 
$$ 
 If $x+D_x z\notin {\rm int}\,K$, then the integrand is reduced to two terms.  Moreover, for $|z|\le 1/\kappa_p$ we have the bound   
$$
|D_x  z pu'(px))|\le |z||p||x|u'(px)\le u(px) 
$$
implying that   
$$
-f(x)-D_x  zf'(x))=-u(px)-D_x z  pu'(px))\le 0. 
$$
We obtain from here, taking into account that $u(px)\ge 0$, the bound 
$$
\cI(f,x)\le u'(px)|p||x|\int_{\{|z|>1/\kappa_p\}}  I_{{\rm int}\, K}(x+D_xz)  \Pi(dz)\le \eta_pu(px).  
$$

Suppose 
that $\langle A(x)p,p\rangle\neq 0$.
Isolating the  full
square we obtain that 
 \bean 
\cL_0 f(x)&=&\frac 12 \left [\langle A(x)p,p\rangle u''(px) + 2\langle
\mu (x),p\rangle u'(px)+ \frac{\langle \mu (x),p\rangle^2}{\langle
A(x)p,p\rangle}\frac {u'^2(px)}{u''(px)} \right]\\ &&
 + \frac 12 \frac{\langle \mu (x),p\rangle^2}{\langle
A(x)p,p\rangle}R(px)u(px) +\cI(f,x)-\beta u(px)\\
&\le&\frac 12 \frac{\langle \mu (x),p\rangle^2}{\langle
A(x)p,p\rangle}R(px)u(px) +\eta_pu(px)-\beta u(px). 
 \eean
It follows that  $ \cL_0 f(x)\le 0$ if $\beta \ge \tilde \eta_p\bar R+\eta_p$. 

Of course, if $\langle A(x)p,p\rangle=0$ we cannot argue as above. 
In this case 
$$
\cL_0 f(x)\le  \langle
\mu (x),p\rangle u'(px)+\eta_p u(px)-\beta u(px). 
$$
Taking into account that 
$$
\sup_{x\in {\rm int}\, K}\frac{\langle \mu (x),p\rangle u'(px)}{u(px)}\le \max_i|\mu_i| \sup_{x\in {\rm int}\, K}\frac {u'(px)}{u(px)}|p||x|= \max_i|\mu_i|\kappa_p, 
$$
we get that $ \cL_0 f(x)\le 0$, if 
$$
\beta\ge \eta_p + \max_i|\mu_i|\kappa_p.  
$$
Combining these two cases, we get the result. \fdem

\smallskip
\noindent
{\bf Remark.} An inspection of the above arguments shows that one can get 
that the  $f_p$ is a Lyapunov  function for 
$$
\beta\ge \sup_{x\in {\rm int}\,K}\left\{\frac 12 \frac{\langle \mu (x),p\rangle^2}{\langle
A(x)p,p\rangle}R(px)I_{\{\langle A(x)p,p\rangle\neq 0\}}+ \frac{\langle \mu (x),p\rangle u'(px)}{u(px)}I_{\{\langle A(x)p,p\rangle\neq 0\}}\right\}+\eta_p. 
$$
Of course, such a bound is less tractable than that given above. 

\smallskip
\noindent
{\bf Construction of classical supersolutions.}
\smallskip
 
Similar argments are useful  in the search of
classical supersolutions for the operator $\cL$. Since $\cL
f=\cL_0 f+U^*(f')$, it is natural to choose $u$ related to $U$.
For a particular case, where $\cC={\bf R}^d_+$ and $U(c)=u(e_1c)$,
with $u$ satisfying the postulated properties (except, maybe,
unboundedness) and assuming, moreover, that the inequality
\beq
\label{cond2} u^*(au'(z))\le g(a)u(z)
\eeq
holds, we get, using
the homogeneity of $\cL_0$, the following result.

 \begin{prop}
\label{supersol}
Let $p\in {\rm int}\,K$. Suppose that
 (\ref{cond2}) holds for every
$a,z>0$ with  $g(a)=o(a)$ as $a\to \infty$. If $\kappa_p<\infty$ and
$\beta< \tilde \eta_p \bar R +\eta_p+\max_i |\mu_i|\kappa_p$, then there exists $a_0$ such that for every $a\ge
a_0$ the function $af_p$ is a  classical strict supersolution of
(\ref{vi}). 
\end{prop}


For the power utility function  $u(z)=z^\gamma/\gamma$, $\gamma\in
]0,1[$, we have:
$$
u^*(au'(z))=(1-\gamma) a^{\gamma/(\gamma-1)}u(z).
$$ 
Therefore, the inequality (\ref{cond2}) holds with $g(a)=o(a)$, $a\to 0$.

\acknowledgement{The research is funded by the grant 14.А12.31.0007. }


\begin{thebibliography}{100}
\bibitem{AMS}
Akian M., Menaldi J.L., Sulem A. On an investment-consumption
model with transaction costs. {\em SIAM J. Control and
Optimization}, {\bf 34} (1996), 1, 329--364.

\bibitem{Alv-Tourin}
Alvarez O., Tourin A.
Viscosity solutions of nonlinear integro-differential equations. {\em Annales de l'institut Henri Poincar\'e. Analyse non lin\'eaire}, {\bf 13} (1996), 3, 293--317.


\bibitem{Arisawa}
Arisawa M.
A new definition of viscosity solutions for a class of second-order degenerate elliptic integro-differential equations.
{\em Annales de l'institut Henri Poincar\'e. Analyse non lin\'eaire}, {\bf 23} (2006), 5, 695--711.

\bibitem{Arisawa2}
 Arisawa M.
A remark on the definitions of viscosity solutions for the integro-differential equations with L\'evy operators. {\em J. Maths. Pures et Applique\'es}, {\bf 89} (2008),  567--574.

\bibitem{Aubin}
Aubin~J.-P.
\newblock {\em Optima and Equilibria. An Introduction to Nonlinear Analysis}.
\newblock Springer,
Berlin--Heidelberg--New York, 1993.

\bibitem{Bardi-Dol}
Bardi M., Dolcetta C. {\sl Optimal Control and Viscosity Solutions
of Hamilton--Jacobi--Bellman Equations}. Birkh\"auser, 1997.

\bibitem{Barles-Soner}
Barles G., Soner H.M. Option pricing with transaction costs and a
non-linear Black-Scholes equation. {\em Finance and Stochastics},
{\bf 2} (1998), 4, 369--397.

\bibitem{Barles-Chasseigne-Imbert}
Barles G., Chasseigne E., Imbert C.
The Dirichlet problem for second-order elliptic integro-differential equations.
{\em Indiana Univ. Math. J.}, {\bf 57} (2008), 1, 213--146.

\bibitem{Barles-Imbert}
Barles G., Imbert C. Second-order elliptic integro-differential equations: viscosity solutions' theory revisited.
{\em Annales de l'institut Henri Poincar\'e. Analyse non lin\'eaire}, {\bf 25} (2008), 3, 567--585.

\bibitem{Benthal}
Benth F.E., Karlsen K.H., Reikvam K. Portfolio optimization in a
L\'evy market with intertemporal substitution
 and transaction costs.
{\sl Stochastics and Stochastics Reports}, {\bf 74} (2002), 3--4, 517--569.

\bibitem{Bentha}
Benth F.E., Karlsen K.H., Reikvam K. Optimal portfolio management rules in a non-Gaussian market with durability and intertemporal substitution.  {\sl Finance and Stochastics}, {\bf 5} (2001), 4, 447--467.

\bibitem{guide}
Crandall M., Ishii H., Lions P.-L. User's guide to viscosity
solutions of second order partial differential equations. {\em
Bulletin of American Mathematical Society}, {\bf 277} (1983),
1-42.

\bibitem{DN}
Davis M., Norman A. Portfolio selection with   transaction costs.
{\em Math. Oper. Res.}, {\bf 15} (1990), 676--713.

\bibitem{FOS}
Framstad N.C., Oksendal B., Sulem A. Optimal consumption and portfolio in a jump diffusion market with proportional transaction costs. {\sl Journal of Mathematical Economics}, {\bf 35} (2001), 233--257.

\bibitem{JK}
Jacobsen E.R., Karlsen K.H.
Continuous dependence estimates for viscosity
solutions of integro-PDEs. J. Differential Equations, {\bf 212} (2005), 278-ñ318.

\bibitem{JS} Jacod J., Shiryaev A.N. Limit theorems for stochastic processes. 2nd edition, Springer, 2002.


\bibitem{Jan-Sh}
Jane\u{c}ek K., Shreve S. Asymptotic analysis for optimal
investment and consumption with transaction costs.
 {\em Finance and Stochastics},  {\bf 8}  (2004),  2, 181--206.

\bibitem{K-Kl}
Kabanov Yu.M., Kl\"uppelberg C. A geometric approach to portfolio
optimization  in models
 with transaction costs. {\sl Finance and
Stochastics}, {\bf 8} (2004), 2, 207--227.

\bibitem{K-Saf} Kabanov Yu.M., Safarian M. {\sl Markets With Transaction Costs: Mathematical Theory.} Springer, Berlin--Heidelberg--New York, 2009.

\bibitem{L-Sh-2}
Liptser R. Sh., Shiryaev A. N. {\sl Theory of Martingales.}
Kluwer, Dordrecht, 1989.

\bibitem{Merton}
Merton R.C. Optimum consumption and portfolio rules in a
continuous time model. {\sl J. Econ. Theory},  {\bf 3} (1971),
373-413.

\bibitem{Pham}
Pham H. Optimal stopping of controlled jump-diffusion processes : a viscosity solutions approach.  {\em Journal of Mathematical Systems, Estimation and Control}, {\bf 8} (1998), 1,

\bibitem{Sayah}
Sayah A. Equations d'Hamilton--Jacobi du premier ordre avec termes int\'egro-diff\'erentiels,  I , II. Comm. P.D.E., {\bf 16} (1991),  1057--1093.

\bibitem{Shreve}
Shreve S. Liquidity premium for capital asset pricing with
transaction costs. In: {\sl Mathematical Finance IMA Volume 65}.
Eds. M.H.A. Davis et al., Springer, New York, 1995, 117--133.

\bibitem{Sh-Son}
Shreve S., Soner M. Optimal investment and consumption  with
transaction costs. {\em Annals Appl. Probab.} {\bf 4} (1994), 3,
609--692.

\bibitem{Soner86a}
Soner H.M.  Optimal control with state space constraint II, {\em SIAM J. Control Optim.}, {\bf 24} (1986), 1110--1122.

\bibitem{Soner86b}
Soner H.M.  Optimal control of jump-markov processes and viscosity solutions. {\em IMA Vols.
in Math. and Applic.}, 10, 501-511.  Springer-Verlag, New York, 1986

\bibitem{TS} Touzi N., Soner  H.M. Dynamic programming for stochastic target problems and geometric flows. {\sl Journal of the European Mathematical Society}, 4, 201-236, 2002.




\end{thebibliography}
\end{document}